\documentclass[twoside, 11pt]{amsart}

\usepackage{amsfonts}
\usepackage[utf8]{inputenc}
\usepackage{array}

\usepackage[english]{babel}
\usepackage{latexsym}
\usepackage{amscd}
\usepackage{amsmath}
\usepackage{amssymb}
\usepackage{amsthm}
\usepackage{graphicx}
\usepackage{mathrsfs}
\usepackage{mathtools}
\usepackage{pict2e}
\usepackage{comment}
\usepackage{stackrel} 
\usepackage{wasysym}
\usepackage[all]{xy}
\usepackage[dvipsnames]{xcolor}
\usepackage[colorlinks,final,hyperindex]{hyperref}
\usepackage{enumitem}
\usepackage{tikz}
\usepackage{tikz-cd}
\usetikzlibrary{decorations.pathmorphing,decorations.markings,arrows,calc,shapes.geometric,arrows.meta,positioning}

\usepackage{geometry}

\allowdisplaybreaks
\usepackage[noabbrev,capitalize]{cleveref}

\newtheorem{lemma}{Lemma}[section]
\newtheorem{theorem}[lemma]{Theorem}

\newtheorem{corollary}[lemma]{Corollary}
\newtheorem{proposition}[lemma]{Proposition}

\newtheorem*{introtheorem}{Theorem}
\theoremstyle{definition}
\newtheorem{definition}[lemma]{Definition}
\newtheorem{remark}[lemma]{\sc Remark}

\renewcommand{\HF}{\mathrm{H} \mathfrak F}
\newcommand{\nocontentsline}[3]{}
\newcommand{\tocless}[2]{\bgroup\let\addcontentsline=\nocontentsline#1{#2}\egroup}

\hypersetup{citecolor=BleuMinuit, linkcolor=bordeau, filecolor=black, urlcolor=BleuMinuit}

\title{Non-formality of prefactorization algebras}

\author{Coline Emprin}
\address{\noindent Coline Emprin, Department of Mathematics, Stockholm University, SE-106 91, Sweden}
\email{\noindent \href{email:coline.emprin@math.su.se}{coline.emprin@math.su.se}}

\author{Alex Takeda}
\address{\noindent Alex Takeda, Instituto de Matem\'atica Pura e Aplicada, 22460-320, Rio de Janeiro, Brazil}
\email{\href{mailto:alex.takeda@impa.br}{alex.takeda@impa.br}}
\date{\today}
\subjclass[2020]{18M70, 18N40, 17B55, 81T45, 81T70}

\begin{document}
		\definecolor{red}{RGB}{230,97,0}
\definecolor{blue}{RGB}{93,58,155}
\definecolor{Chocolat}{rgb}{0.36, 0.2, 0.09}
\definecolor{BleuTresFonce}{rgb}{0.215, 0.215, 0.36}
\definecolor{BleuMinuit}{RGB}{0, 51, 102}
\definecolor{bordeau}{rgb}{0.5,0,0}
\definecolor{BORDEAU}{rgb}{0.5,0,0}
\definecolor{turquoise}{RGB}{6, 62, 62}
\definecolor{rose}{RGB}{235, 62, 124}

\newcommand{\coline}[1]{\textcolor{rose}{#1}}

\newcommand{\hooklongrightarrow}{\lhook\joinrel\longrightarrow}

\newcommand{\chdiscr}{\mathsf{discr.\ Ch}}
\newcommand{\sSe}{\mathsf{sSet}}
\newcommand{\sLialg}{\ensuremath{\mathcal{sL}_\infty\text{-}\,\mathsf{alg}}}
\newcommand{\isLialg}{\ensuremath{\infty\text{-}\,\mathcal{sL}_\infty\text{-}\,\mathsf{alg}}}
\newcommand{\sLiealg}{\ensuremath{\mathcal{s}\lie\,\text{-}\,\mathsf{alg}}}
\def\cD{\Delta}

\newcommand{\ucomnalg}{\ensuremath{\mathrm{uCom}_{\leslant 0}\text{-}\,\mathsf{alg}}}
\newcommand{\fnQsp}{\ensuremath{\mathsf{ft}\mathbb{Q}\text{-}\mathsf{ho}(\mathsf{sSet})}}
\newcommand{\fnsp}{\ensuremath{\mathsf{fn}\text{-}\mathsf{Sp}}}
\newcommand{\nQsp}{\ensuremath{\mathsf{n}\mathbb{Q}\text{-}\mathsf{Sp}}}
\newcommand{\nsp}{\ensuremath{\mathsf{n}\text{-}\mathsf{Sp}}}
\newcommand{\fQalg}{\ensuremath{\mathsf{ft}\text{-}\mathsf{ho}(\mathrm{uCom}_{\leqslant 0}\text{-}\,\mathsf{alg})}}

\newcommand{\coker}{\operatorname{coker}}
\newcommand{\Ran}{\operatorname{Ran}}
\newcommand{\ho}{\operatorname{ho}}
\newcommand{\Mal}{\operatorname{Mal}}
\newcommand{\Loc}{\operatorname{Loc}}

\newcommand{\sLi}{\ensuremath{\mathcal{sL}_\infty}}
\newcommand{\ccoLi}{\ensuremath{\mathrm{cBCom}}}
\newcommand{\sLie}{\ensuremath{\mathcal{s}\mathrm{Lie}}}
\newcommand{\com}{\ensuremath{\mathrm{Com}}}
\newcommand{\ucom}{\ensuremath{\mathrm{uCom}}}
\newcommand{\Cobar}{\ensuremath{\Omega}}
\newcommand{\hatCobar}{\ensuremath{\widehat{\Omega}}}
\renewcommand{\Bar}{\ensuremath{\mathrm{B}}}
\newcommand{\hatBar}{\ensuremath{\widehat{\mathrm{B}}}}
\newcommand{\lie}{\ensuremath{\mathrm{Lie}}}
\newcommand{\uhocom}{\ensuremath{\mathrm{u}\Omega\mathrm{BCom}}}

\newcommand{\Sh}{\ensuremath{\mathrm{Sh}}}

\newcommand{\g}{\ensuremath{\mathfrak{g}}}
\newcommand{\wsLi}{\ensuremath{\widehat{\mathcal{sL}_\infty}}}

\newcommand{\Rh}{\ensuremath{\mathrm{R^h}}}
\renewcommand{\L}{\ensuremath{\mathscr{L}}}
\newcommand{\Li}{\mathfrak{L}}
\newcommand{\APL}{\mathrm{A_{PL}}}
\newcommand{\CPL}{\mathrm{C_{PL}}}

\newcommand{\antishriek}{\text{\raisebox{\depth}{\textexclamdown}}}
\newcommand{\RT}{\mathrm{RT}}
\newcommand{\LRT}{\mathrm{LRT}}
\newcommand{\Sy}{\mathbb{S}}
\renewcommand{\d}{\ensuremath{\mathrm{d}}}
\newcommand{\PP}{\ensuremath{\mathrm{P}}}
\def\Ho#1#2{\Lambda^{#2}_{#1}}
\def\De#1{\Delta^{#1}}
\newcommand{\NN}{\mathbb{N}}
\newcommand{\RR}{\mathbb{R}}
\def\BCH{\mathrm{BCH}}
\newcommand{\wPT}{\ensuremath{\mathrm{wPT}}}
\newcommand{\oPT}{\ensuremath{\overline{\mathrm{PT}}}}
\newcommand{\PaRT}{\ensuremath{\mathrm{PaRT}}}
\newcommand{\PaPT}{\ensuremath{\mathrm{PaPT}}}
\newcommand{\PaPRT}{\ensuremath{\mathrm{PaPRT}}}
\def\rmC{\mathrm{C}}
\newcommand{\berglund}{\ensuremath{\mathcal{B}}}
\def\hot{\widehat{\otimes}} 
\def\whk{\widehat{k}}

\def\colim{\mathop{\mathrm{colim}}}

\newcommand{\CC}{\ensuremath{\mathrm{CC}_\infty}}

\newcommand{\Lalg}{\ensuremath{\mathscr{L}_\infty\text{-}\mathsf{alg}}}

\newcommand{\F}{\ensuremath{\mathrm{F}}}
\newcommand{\h}{\ensuremath{\mathfrak{h}}}

\newcommand{\QQ}{\ensuremath{\mathbb{Q}}}
\newcommand{\D}{\ensuremath{\mathscr{D}}}
\renewcommand{\P}{\ensuremath{\mathcal{P}}}
\newcommand{\C}{\ensuremath{\mathcal{C}}}
\newcommand{\VdL}{\ensuremath{\mathrm{VdL}}}
\newcommand{\ch}{\ensuremath{\mathrm{Ch}}}
\newcommand{\End}{\ensuremath{\mathrm{End}}}
\newcommand{\Aut}{\ensuremath{\mathrm{Aut}}}
\newcommand{\eend}{\ensuremath{\mathrm{end}}}
\newcommand{\susp}{\ensuremath{\mathscr{S}}}
\newcommand{\T}{\ensuremath{\mathcal{T}}}
\newcommand{\vdl}{\ensuremath{\mathrm{VdL}}}
\newcommand{\Tw}{\ensuremath{\mathrm{Tw}}}
\newcommand{\Hom}{\ensuremath{\mathrm{Hom}}}
\renewcommand{\S}{\ensuremath{\mathbb{S}}}
\renewcommand{\k}{\ensuremath{\mathbb{K}}}
\newcommand{\id}{\ensuremath{\mathrm{id}}}
\newcommand{\MC}{\ensuremath{\mathrm{MC}}}
\newcommand{\mc}{\ensuremath{\mathfrak{mc}}}
\newcommand{\mclie}{\ensuremath{\overline{\mathfrak{mc}}}}
\newcommand{\dgl}{\ensuremath{\mathsf{dgLie}}}

\newcommand{\B}{\mathcal{B}}
\newcommand{\Z}{\mathbb{Z}}

\newcommand{\ad}{\operatorname{ad}}
\newcommand{\PTt}{\ensuremath{\widetilde{\mathrm{PT}}}}
\newcommand{\PT}{\ensuremath{\mathrm{PT}}}

\newcommand{\A}{\ensuremath{\mathrm{A}}}

\makeatletter

\newcommand{\R}{\mathbb{R}}
\newcommand{\K}{\mathcal{K}}
\newcommand{\N}{\mathbb{N}}
\newcommand{\calA}{\mathcal{A}}
\newcommand{\calB}{\mathcal{B}}
\newcommand{\calC}{\mathcal{C}}
\newcommand{\calD}{\mathcal{D}}
\newcommand{\calE}{\mathcal{E}}
\newcommand{\calF}{\mathcal{F}}
\newcommand{\calG}{\mathcal{G}}
\newcommand{\calH}{\mathcal{H}}
\newcommand{\calI}{\mathcal{I}}
\newcommand{\calJ}{\mathcal{J}}
\newcommand{\calK}{\mathcal{K}}
\newcommand{\calL}{\mathcal{L}}
\newcommand{\calM}{\mathcal{M}}
\newcommand{\calN}{\mathcal{N}}
\newcommand{\calO}{\mathcal{O}}
\newcommand{\calP}{\mathcal{P}}
\newcommand{\calQ}{\mathcal{Q}}
\newcommand{\calR}{\mathcal{R}}
\newcommand{\calS}{\mathcal{S}}
\newcommand{\calT}{\mathcal{T}}
\newcommand{\calU}{\mathcal{U}}
\newcommand{\calV}{\mathcal{V}}
\newcommand{\calW}{\mathcal{W}}
\newcommand{\calX}{\mathcal{X}}
\newcommand{\calY}{\mathcal{Y}}
\newcommand{\calZ}{\mathcal{Z}}
\newcommand{\kk}{{\Bbbk}}

\newcommand{\ev}{\mathop{\mathrm{ev}}\nolimits}
\newcommand{\co}{\mathop{\mathrm{co}}\nolimits}
\newcommand{\pt}{\mathop{\mathrm{pt}}\nolimits}
\newcommand{\wt}{\mathop{\mathrm{wt}}\nolimits}
\newcommand{\pCY}{\mathop{\mathrm{pCY}}\nolimits}
\newcommand{\gr}{\mathop{\mathrm{gr}}\nolimits}

\newcommand{\Imm}{\mathop{\mathrm{Im}}\nolimits}

\tikzset{>={Stealth[scale=1.2]}}
\tikzset{->-/.style={decoration={
			markings,
			mark=at position #1 with {\arrow{>}}},postaction={decorate}}}
\tikzset{-w-/.style={decoration={
			markings,
			mark=at position #1 with {\arrow{Stealth[fill=white,scale=1.4]}}},postaction={decorate}}}
\tikzset{->-/.default=0.65}
\tikzset{-w-/.default=0.65}
\tikzstyle{bullet}=[circle,fill=black,inner sep=0.5mm]
\tikzstyle{circ}=[circle,draw=black,fill=white,inner sep=0.5mm]
\tikzstyle{vertex}=[circle,draw=black,thick,inner sep=0.5mm]
\tikzstyle{dot}=[draw,circle,fill=black,minimum size=0.5mm,inner sep = 0mm, outer sep = 0mm]
\tikzset{darrow/.style={double distance = 4pt,>={Implies},->},
	darrowthin/.style={double equal sign distance,>={Implies},->},
	tarrow/.style={-,preaction={draw,darrow}},
	qarrow/.style={preaction={draw,darrow,shorten >=0pt},shorten >=1pt,-,double,double
		distance=0.2pt}}
\newcommand\tikcirc[1][2.5]{\tikz[baseline=-#1]{\draw[thick](0,0)circle[radius=#1mm];}}
\newcommand{\tikzfig}[1]{\begin{tikzpicture}[auto,baseline={([yshift=-.5ex]current bounding box.center)}]#1\end{tikzpicture}}
\usetikzlibrary{decorations.pathreplacing}

\newcommand{\DD}{\EuScript D}
\newcommand{\BB}{\EuScript B}
\newcommand{\OO}{\EuScript O}
\newcommand{\KK}{\EuScript K}
\newcommand{\YY}{\EuScript Y}
\newcommand\nn         {\nonumber \\}
\newcommand{\pf}{{\it Proof.}\hspace{2ex}}
\newcommand{\apf}{{\it Another Proof.}\hspace{2ex}}
\newcommand{\epf}{\hspace*{\fill}\emph{$\square$}}
\newcommand{\halmos}{\rule{1ex}{1.4ex}}
\newcommand{\pfbox}{\hspace*{\fill}\emph{$\halmos$}}

\newcommand{\ac}{\mathop{ac}\nolimits}
\newcommand{\sC}{\mathop{s^{-1}\overline{C}}\nolimits}
\newcommand{\sA}{\mathop{s\overline{A}}\nolimits}
\newcommand{\Rep}{\mathop{\mathrm{Rep}}\nolimits}
\newcommand{\Span}{\mathop{\mathrm{Span}}\nolimits}
\newcommand{\proj}{\mathop{\mathrm{proj}}\nolimits}
\newcommand{\PHom}{\mathop{\mathrm{PHom}}\nolimits}
\newcommand{\Dim}{\mathop{\mathrm{Dim}}\nolimits}

\newcommand{\Aus}{\mathop{\mathrm{Aus}}\nolimits}
\newcommand{\Ann}{\mathop{\mathrm{Ann}}\nolimits}
\newcommand{\APC}{\mathop{\mathrm{APC}}\nolimits}
\newcommand{\Barr}{\mathop{\mathrm{Bar}}\nolimits}
\newcommand{\Cone}{\mathop{\mathrm{Cone}}\nolimits}
\newcommand{\modu}{\mathop{\mathrm{mod}}\nolimits}
\newcommand{\DGMod}{\mathop{\mathrm{DGMod}}\nolimits}
\newcommand{\soc}{\mathop{\mathrm{soc}}\nolimits}
\newcommand{\dg}{\mathop{\mathrm{dg}}\nolimits}
\newcommand{\red}{\mathop{\mathrm{red}}\nolimits}
\newcommand{\Map}{\mathop{\mathrm{Map}}\nolimits}
\newcommand{\inj}{\mathop{\mathrm{inj}}\nolimits}
\newcommand\M          {{\mathcal{M}}}
\newcommand{\Ker}{\mathop{\mathrm{Ker}}\nolimits}
\newcommand{\Tor}{\mathop{\mathrm{Tor}}\nolimits}
\newcommand{\Tot}{\mathop{\mathrm{Tot}}\nolimits}
\newcommand{\Ext}{\mathop{\mathrm{Ext}}\nolimits}
\newcommand{\sg}{\mathop{\mathrm{sg}}\nolimits}
\newcommand{\Sl}{\mathop{\mathfrak{sl}}\nolimits}
\newcommand{\nc}{\mathop{\mathrm{nc}}\nolimits}
\newcommand{\Tr}{\mathop{\mathrm{Tr}}\nolimits}
\newcommand{\Irr}{\mathop{\mathrm{Irr}}\nolimits}
\newcommand{\HC}{\mathop{\mathrm{HC}}\nolimits}
\newcommand{\HH}{\mathop{\mathrm{HH}}\nolimits}
\newcommand{\THH}{\mathop{\mathrm{TH}}\nolimits}
\newcommand{\Perf}{\mathop{\mathrm{Perf}}\nolimits}
\newcommand{\del}{\partial}

\newcommand{\verteq}{\rotatebox{90}{$\;\;=\;\;$}}

\newcommand{\Sym}{\operatorname{Sym}}
\newcommand{\CE}{\operatorname{CE}}
\newcommand{\Disk}{\operatorname{Disk}}
\newcommand{\Ch}{\operatorname{Ch}}
\newcommand{\Alg}{\operatorname{Alg}}
\newcommand{\dR}{\mathrm{dR}}
\newcommand{\cpt}{\mathrm c}
\newcommand{\Pfac}{\mathcal P}
\newcommand{\PF}{\mathfrak F}
\newcommand{\Ug}{\mathcal U\mathfrak g}
\newcommand{\op}[2]{\iota^{#1}_{#2}}
\newcommand{\bbD}{\underline D}
\newcommand{\bbE}{\mathbf E}
\newcommand{\tensor}{\mathbin{\otimes}}
\newcommand{\xto}[1]{\xrightarrow{#1}}
\newcommand{\xfrom}[1]{\xleftarrow{#1}}
\newcommand{\simeqto}{\xrightarrow{\sim}}
\newcommand{\bracket}[2]{\left\{#1,#2\right\}}
\newcommand{\tb}{\widetilde\beta}

\begin{abstract}
We construct a Poisson bracket on the cohomology of locally constant
prefactorization algebras on \(\mathbb{R}^m\) that provides an obstruction
to formality in every dimension \(m\geqslant 2\). Our construction combines
colored obstruction sequences to formality with linking spheres and
generalizes the two-dimensional bracket of Bruinsma, Schenkel, and Vicedo.
We use it to extend their non-formality theorem for Lie algebra
factorization envelopes to all higher dimensions. We also establish new
non-formality results for factorization envelopes of Lie algebras, their compactifications and generalizations to \(L_\infty\)-algebras, free BF-type theories, higher abelian
Chern--Simons theories and their compactifications.
\end{abstract}

	\maketitle

\makeatletter
\def\@tocline#1#2#3#4#5#6#7{\relax
	\ifnum #1>\c@tocdepth %
	\else
	\par \addpenalty\@secpenalty\addvspace{#2}%
	\begingroup \hyphenpenalty\@M
	\@ifempty{#4}{%
		\@tempdima\csname r@tocindent\number#1\endcsname\relax
	}{%
		\@tempdima#4\relax
	}%
	\parindent\z@ \leftskip#3\relax \advance\leftskip\@tempdima\relax
	\rightskip\@pnumwidth plus4em \parfillskip-\@pnumwidth
	#5\leavevmode\hskip-\@tempdima
	\ifcase #1
	\or\or \hskip 1em \or \hskip 2em \else \hskip 3em \fi%
	#6\nobreak\relax
	\hfill\hbox to\@pnumwidth{\@tocpagenum{#7}}\par%
	\nobreak
	\endgroup
	\fi}

\newcommand{\enableopenany}{%
	\@openrightfalse%
}
\makeatother
	
	\setcounter{tocdepth}{1}
	\tableofcontents

\section{\textcolor{bordeau}{Introduction}}
\noindent Prefactorization algebras are certain algebraic objects associated to a manifold that are proposed to give a framework for local observables of quantum field theories \cite{CG17,CG21}. These objects axiomatize the idea that a quantum field theory assigns observables to each region in spacetime, and inclusions of disjoint regions give rise to operations between the corresponding spaces of observables. In the setting considered here, these regions are smoothly embedded open disks, and a prefactorization algebra $\mathfrak F$ assigns a cochain complex $\mathfrak F(D)$ to each open disk $D$, together with compatible operations combining observables supported in disjoint disks into an observable supported in a larger disk. Spaces of observables given by cochain complexes arise naturally in the Batalin--Vilkovisky formalism for quantization of field theories, where the differential can account for the presence of gauge symmetries. 

\medskip

\noindent Equivalently, we consider algebras in cochain complexes over a certain colored operad $\Pfac_M$ encoding configurations of disks embedded in $M$. Taking cohomology of such an algebra yields a prefactorization algebra $\HF$ valued in graded vector spaces. In general, this cohomology prefactorization algebra $\HF$ does not capture the full homotopical information carried by $\mathfrak F$. Roughly speaking, it only remembers how two observables are combined by a simple product, and not how they are related as their supports move with respect to one another.

\medskip

\noindent The prefactorization algebra $\mathfrak F$ is called \emph{formal} if it is connected to $\HF$, equipped with its induced $\Pfac_M$-algebra structure, by a zigzag of quasi-isomorphisms of prefactorization algebras; in these special cases, $\HF$ does remember all the homotopical information of $\mathfrak F$. Applying to this setting the formalism of Dimitrova \cite{Dim12}, Bruinsma, Schenkel and Vicedo obtain a first obstruction class $\vartheta_2$ to the formality of a prefactorization algebra \cite{BSV24}. A case of particular interest is that of \emph{locally constant} prefactorization algebras,  where every inclusion of disks induces a quasi-isomorphism. Consequently, its cohomology assigns to every disk a value isomorphic to \(H=\mathcal{H}\mathfrak{F}(\mathbb{R}^m)\). In this case, when the ambient dimension is $m=2$, they extract from $\vartheta_2$ a map of degree $-1$
\[ L^{(2)} \colon H \otimes H \to H \ , \]
which makes $H$ into a $(-1)$-shifted Poisson algebra. The bracket  can be geometrically interpreted as a linking operation in two dimensions, where one observable is moved along a linking circle to another. They prove that non-vanishing of the map $L^{(2)}$ implies non-vanishing of the class $\vartheta_2$, and therefore non-formality of $\mathfrak F$. 

\medskip

\noindent Our goal in this paper is to explore higher dimensional obstruction classes for the formality of prefactorization algebras, having in mind the following organizing question:   \textit{How much more data does a topological theory (in its formulation as a locally constant prefactorization algebra) carry, besides the cohomology product and the linking operation?}   

\medskip

\noindent The example  we examine in \cref{prop:massey-counterexample} shows that the linking invariant $L^{(2)}$ does not capture the entire obstruction. In this example, starting from a dg algebra that has a nontrivial triple Massey product, we give a locally constant prefactorization algebra in any dimension that always has $\vartheta_2 \neq 0$ in any spatial dimension, but where $L^{(2)}=0$ in dimension two.

\medskip

\noindent The first step towards a systematic answer to our motivating question is to understand what are the higher obstruction classes for formality of a general prefactorization algebra, even without assuming local constancy. Building on the work of Kaledin~\cite{Kal07}, Lunts~\cite{Lun07}, and
Melani and Rubi\'o~\cite{MR19}, the first-named author extended the
Kaledin class to algebras over colored operads~\cite{Kaledin}. This
construction associates to each prefactorization algebra $\mathcal{F}$ a single
cohomology class, its Kaledin class, which vanishes if and only if $\mathcal{F}$
is formal. Its very first truncation recovers Dimitrova's class
$\vartheta_2$. Here we study the
higher-order obstructions through successive truncations of the Kaledin class, following the inductive approach of~\cite{CE24b}. If
$\vartheta_m$ vanishes, one can define the next obstruction to
formality, $\vartheta_{m+1}$.
\medskip

\noindent  For locally constant prefactorization algebras, the higher
obstruction classes $\vartheta_m$ allow us to extend the bracket
construction of \cite{BSV24} to every
dimension. For all $m\geqslant 2$, we replace the loop by an $(m-1)$-sphere surrounding a disk, described by chains of disk inclusions. Using a more systematic description of the relation between posets of disks and homologies of configuration spaces, we obtain the following higher-dimensional generalization.

\medskip 

\begin{introtheorem}[Proposition~\ref{prop:higher-poisson-invariant} and Theorem~\ref{thm:higher-poisson-bracket}]
Let $\mathfrak F$ be a locally constant prefactorization algebra on $\R^m$, with $m\geqslant2$, such that $\vartheta_2 = \cdots = \vartheta_{m-1} = 0$. 
\begin{enumerate}
    \item $ \HF(\R^m)$ carries a Poisson bracket $L^{(m)}$ of degree $1-m$  determined by $\vartheta_m$; 
    \item  if $L^{(m)} $ is nonzero, then $\mathfrak F$ is not formal.
\end{enumerate}
\end{introtheorem}

\medskip

\noindent We compute this bracket for several families of examples, and use these computations to establish the non-formality of a wide class of prefactorization algebras coming from free constructions.

\bigskip

\paragraph{\bf Factorization envelopes}
The factorization envelope of a local Lie algebra is constructed from the Chevalley--Eilenberg chains of its compactly supported sections \cite[Section~3.6]{CG17}. For a finite-dimensional Lie algebra $\mathfrak g$ over $\R$ or $\mathbb{C}$, its factorization envelope is given by the symmetric algebra over 
\[\mathfrak g\otimes\Omega^\bullet_{\R^m} \] 
endowed with the usual symmetric product. We prove in Theorem~\ref{thm:factorization-envelope-higher-bracket} that its bracket $L^{(m)}$ is given by $(-1)^m$ times the shifted Poisson bracket induced by the Lie bracket of $\mathfrak g$, and hence non-zero when $\mathfrak g$ is non-abelian. We extend this construction to uncurved $\mathcal{L}_\infty$-algebra, where an analogous calculation holds; we summarize our results for these cases in the following theorem. 

\begin{introtheorem}[Corollary~\ref{cor:factorization-envelope-equivalence}, Theorem~\ref{thm:coefficient-envelope-first-invariant} and Proposition~\ref{prop:compactified-abelian-formal}] Let $m\geqslant2$.
    \begin{enumerate}[label=\textup{(\roman*)}]
        \item For a finite-dimensional real or complex Lie algebra $\mathfrak g$ in degree zero, its factorization envelope on $\R^m$ is formal if and only if $\mathfrak g$ is abelian. The same criterion holds for the pushforward from $\R^m\times N$ to $\R^m$ along any nonempty closed manifold $N$.
        \item For an uncurved $L_\infty$-algebra $\mathfrak g$ over $\R$ or $\mathbb C$, non-vanishing of the induced Lie bracket on $H^\bullet(\mathfrak g)$ implies non-formality of the factorization envelope of $\mathfrak g$ on $\R^m$.
    \end{enumerate}
\end{introtheorem}

 \medskip

\paragraph{\bf Free theories from symmetric pairings}
A complex $V$ equipped with a closed graded symmetric pairing of the appropriate degree gives rise to a polynomial prefactorization algebra over $\R[\hbar]$, following the construction of free Batalin--Vilkovisky observables \cite[Section~3.1]{GRW21}. Its bracket $L^{(m)}$ is simply determined by the induced pairing on cohomology of $V$, as follows: 

\begin{introtheorem}[Theorem~\ref{thm:paired-complex-first-invariant}]
    Let $V$ be a  cochain complex with a closed  symmetric pairing 
    \[B_V\colon V\otimes V\to\R[m+1] \ , \] where $m\geqslant2$. 
    The associated polynomial model $\mathfrak F_V$ of free Batalin--Vilkovisky observables on $\R^m$, with coefficients in $\R[\hbar]$, is formal as a real prefactorization algebra if and only if  the induced pairing on cohomology is zero.
\end{introtheorem}

\noindent This applies for free BF theory \cite{CR01,CMR18}, higher abelian Chern--Simons theory \cite[Section~5.4]{GRW21}, and their compactifications over nonempty closed oriented manifolds. For linear Chern--Simons theory on $\R^2\times S^1$ with its standard real pairing, we recover the calculation of \cite[Theorem~5.1]{BSV24}.

\bigskip

\paragraph{\bf AI use statement} The structure of this paper and its conceptual ideas were developed without the use of AI tools. We later used it to suggest examples, such as the counterexample in Proposition~\ref{prop:massey-counterexample} and the applications of \mbox{Subsections~\ref{subsec:dg-Linfty-envelopes} and~\ref{subsec:free-BF}}, and to correct errors in our arguments and simplify the proofs. All statements produced by this process have been verified by hand.

\bigskip

\paragraph{\bf Acknowledgments} 
We would like to thank Victor Carmona for fruitful conversations. The first author was supported by the Göran Gustafsson Foundation for Research in Natural Sciences and Medicine through Alexander Berglund's 2024 Göran Gustafsson Prize. The second author was supported by the Knut and Alice Wallenberg Foundation and by IMPA.

\newpage

\section{\textcolor{bordeau}{Obstruction sequences to formality}}
\label{sec:obstruction-sequences}

\noindent In this section, we recall the obstruction sequences to formality of algebras over colored operads, as established in \cite{CE24b}. This theory is applied in \cref{sec:locally-constant} to the prefactorization operad of a smooth manifold $M$.

\subsection{Colored operads and their algebras}
We start by recalling the definition of a colored operad and associated homotopy coherent algebras with the aim of fixing notations and conventions. For further details, we refer the reader to \cite{LodayVallette12,War19, BSV24}. Let us fix a set of colors $C$ and let $\Bbbk$ be a characteristic zero field.

\begin{definition}[Colored operads and algebras]
\label{def:colored-operad}
A $C$-colored dg operad $\mathcal P$ consists of cochain complexes $\mathcal P(c,\underline c)$ indexed by an output color $c$ and an input profile $\underline c=(c_1,\ldots,c_n)$. It is equipped with symmetric actions permuting the inputs and their colors, closed degree-zero units $1_c\in\mathcal P(c,c)$, and partial composition maps
\[
\circ_i\colon
\mathcal P(c,\underline c)\otimes\mathcal P(c_i,\underline d)
\longrightarrow
\mathcal P(c,\underline c\circ_i\underline d) \ ,
\] where  $\underline c\circ_i\underline d$ is the profile obtained by
replacing $c_i$ in $\underline c$ with  $\underline d$. These structure maps commute with the differentials and satisfy the associativity, unit and equivariance axioms. A morphism $f\colon\mathcal P\to\mathcal Q$ is a family of cochain maps compatible with permutations, units and compositions:
\[\mathcal P(c,\underline c)
\to\mathcal Q(c,\underline c) \ . \]  For a family of complexes $A=(A_c)_{c\in C}$, its endomorphism operad is given by
\[
\operatorname{End}_A(c,\underline c)
=\operatorname{Hom}_{\Bbbk}(A_{c_1}\otimes\cdots\otimes A_{c_n},A_c) \ ,
\]
with composition given by substitution. A \emph{$\mathcal P$-algebra structure} on $A$ is a morphism $\alpha \colon\mathcal P\to\operatorname{End}_A$ of $C$-colored dg operads.
\end{definition}

\noindent From now on, let us fix a family $A=(A_c)_{c\in C}$ and a $C$-colored operad $\mathcal P$ concentrated in degree zero and equipped with an augmentation $\epsilon\colon\mathcal P\to I$, where $I$ is the unit operad on $C$, with $I(c,c)=\Bbbk\,1_c$. Write $\overline{\mathcal P}:=\ker\epsilon$ for the augmentation ideal.

\begin{definition}[Bar construction]
\label{def:bar-construction}
The \emph{bar construction} of $\mathcal P$ is the conilpotent dg cooperad
\[
B\mathcal P
=
\bigl(T^c\left(\overline{\mathcal P}[1]\right),d_{B\mathcal P}\bigr) \ .
\]
Here $T^c$ denotes the cofree conilpotent colored cooperad: its elements are represented by rooted trees with numbered inputs and colored edges, whose (internal) vertices are decorated by elements of $\overline{\mathcal P}[1]$ with matching input and output colors. Decorations are taken modulo the symmetric actions at the vertices. We write $t(u_1,\ldots,u_r)$ for a decorated tree as follows

\[
t(u_1,u_2,u_3,u_4)
=
\begin{tikzpicture}[
  baseline=(current bounding box.center),
  every node/.style={inner sep=1.5pt}
]
\node (u1) at (0,0) {$u_1$};
\node (u2) at (-1.35,-1) {$u_2$};
\node (u3) at (1.15,-1) {$u_3$};
\node (u4) at (1.85,-2) {$u_4$};

\draw (0,.65) -- (u1);
\draw (u1) -- (u2);
\draw (u1) -- (u3);
\draw (u3) -- (u4);

\draw (u2) -- (-2,-1.9)
  ;
\draw (u2) -- (-1.35,-1.9)
   ;
\draw (u2) -- (-.7,-1.9)
  ;
\draw (u3) -- (.45,-1.9)
  ;
\draw (u4) -- (1.2,-2.9)
  ;
\draw (u4) -- (1.85,-2.9)
 ;
\draw (u4) -- (2.5,-2.9)
 ;
\end{tikzpicture}
\]

\noindent The \emph{weight} of a decorated tree is
its number of vertices, leading to the weight decomposition: \[
B\mathcal P=I\oplus\bigoplus_{r\geqslant1}(B\mathcal P)^{(r)} \ .
\]
The differential $d_{B\mathcal P}$ is the signed sum of contractions of internal edges, composing the two adjacent decorations using the partial compositions of $\mathcal P$, see \cite[Sections~5.8 and~6.5--6.6]{LodayVallette12} and \cite[Section~2.6]{War19}.
\end{definition}

\begin{definition}[Convolution dg Lie algebra]
The \emph{colored convolution complex}
\[
\mathfrak g_A
:=
\operatorname{Hom}_{C}
\bigl(\overline{B\mathcal P},\operatorname{End}_A\bigr) \ ,
\]  consists of families
of graded maps over all color profiles, equivariant under permutations of the inputs.  Its differential is given by
\[
\partial\alpha
=
d_{\operatorname{End}_A}\circ\alpha
-
(-1)^{|\alpha|}\alpha\circ d_{B\mathcal P}  \ .
\] The pre-Lie product $\alpha\star\beta$ is computed on a tree $T$
by cutting one internal edge at a time. For each cut, apply $\alpha$ to the part containing the root and $\beta$ to the detached subtree, then insert the second operation into the first at the cut input. The signed sum over cuts gives the pre-Lie product \cite[Section~6.4]{LodayVallette12}. The graded commutator
\[
[\alpha,\beta]
=\alpha\star\beta-(-1)^{|\alpha||\beta|}\beta\star\alpha
\]
makes $\mathfrak g_A$ a dg Lie algebra. The vertex weight on $B\mathcal P$ induces a completed weight decomposition of $\mathfrak g_A$, whose component $\alpha^{(r)}$ is the restriction of $\alpha$ to trees with $r$ vertices. The pre-Lie product adds weights, and in the differential, the term induced by $d_A$ preserves weight, while the term induced by $d_{B\mathcal P}$ raises it by one.
\end{definition}

\begin{definition}[$\mathcal P_\infty$-algebra structures]

A \emph{dg $\mathcal P_\infty$-algebra}  $(A, \alpha)$ is given by a Maurer--Cartan element in $\mathfrak g_A$ of degree one, i.e.
\[
\partial\alpha+ \alpha \star \alpha=0 \ .
\] The structure is \emph{strict} if $\alpha^{(r)}=0$ for
$r\geqslant2$, equivalently if it is induced by a $\mathcal P$-algebra structure. For another family $B$, let us write
\[
\operatorname{End}_B^A(c, \underline c )
:=
\operatorname{Hom}_{\Bbbk}
(A_{c_1}\otimes\cdots\otimes A_{c_n},B_c) \ .
\]
An \emph{$\infty$-morphism} $f\colon(A,\alpha)\rightsquigarrow(B,\beta)$ is a degree-zero element in
\[
\mathfrak a^A_B
:=
\operatorname{Hom}_{C}
\bigl(B\mathcal P,\operatorname{End}^A_B\bigr) \ ,
\]
satisfying
\[
\partial f=f\star\alpha-\beta \circledcirc f \ .
\] Here the differential $\partial$ and $\star$ are defined as in $\mathfrak{g}_A$ whereas $\beta \circledcirc f$ is the sum over the decompositions of $T$
into a root subtree and the subtrees grafted onto its inputs. Apply $f$ to each input subtree and insert the resulting operations into the value of $\beta$ on the root subtree, with the Koszul signs. It is an \emph{$\infty$-quasi-isomorphism} if $f^{(0)}$ is a quasi-isomorphism at every color.
\end{definition}

\begin{remark}
An $\infty$-morphism is a strict $\mathcal P$-algebra morphism if $f^{(r)} = 0$ for all $r \geqslant 1$.
\end{remark}

\subsection{Obstruction classes to formality}
Let $(A,\alpha)$ be a $\mathcal P$-algebra, and let $\alpha_*$ denote the induced structure on $H:=H(A)$. Recall that $(A,\alpha)$ is said \emph{formal} if it is connected to $(H,\alpha_*)$ by a zigzag of $\mathcal P$-algebra morphisms that are quasi-isomorphisms at every color. To construct the obstructions to formality, we use homotopy transfer to obtain a minimal $\mathcal P_\infty$-structure on $H$. Its weight-one component recovers $\alpha_*$, and its higher components give rise to successive obstructions. Since $\Bbbk$ is a field, we may choose at every color a \emph{strong deformation retract}
\[
\begin{tikzcd}[column sep=large]
H \arrow[r,shift left=.7ex,"i"] &
A \arrow[l,shift left=.7ex,"p"] \arrow[loop right,"h"]
\end{tikzcd}
\]
where $i$ and $p$ are cochain maps and $h$ has degree $-1$, with
\[
pi=\operatorname{id}_H \ ,\qquad
d_Ah+hd_A=ip-\operatorname{id}_A \ ,\qquad
hi=ph=h^2=0 \ .
\]
The homotopy transfer construction of \cite[Sections~10.3--10.4]{LodayVallette12} applies colorwise, since its tree formulas involve only compositions along matching colors, see \cite[Section~2.8]{War19}.

\begin{theorem}[Homotopy transfer]
\label{thm:homotopy-transfer}
There exists an induced $\mathcal P_\infty$-structure $(H, \beta)$ for which $i$ extends to an $\infty$-quasi-isomorphism
\[
i_\infty\colon(H,\beta)\rightsquigarrow(A,\alpha) \ .
\]
The algebra $(H,\beta)$ is a \emph{minimal model} of $(A,\alpha)$, unique up to $\infty$-isomorphism.
\end{theorem}

\noindent The transferred structure is given by the tree formulas of \cite[Section~10.3]{LodayVallette12}. The weight-one component of $\beta$ recovers the algebra structure induced on cohomology. Indeed, for a corolla decorated by the suspension $sv$, with $v\in\overline{\mathcal P}(c,\underline c)$, we have
\[
\beta^{(1)}(sv)
=
p_c\,\alpha(v)
\bigl(i_{c_1}\otimes\cdots\otimes i_{c_n}\bigr)
=
\alpha_*(v) \ . 
\] The higher components are computed on the same decorated trees.
For $T=t(u_1,\ldots,u_r)$, write $u_j=sv_j$. To evaluate $\beta^{(r)}(T)$, we apply $\alpha(v_j)$ at vertex $j$, place $i$ at the inputs, $h$ on each internal edge and $p$ at the output, and compose with the suspension and Koszul signs. For the above tree, this gives the following operation
\[
\beta\bigl(t(u_1,u_2,u_3,u_4)\bigr)
=
\begin{tikzpicture}[
  baseline=(v2.base),
  line width=.55pt, line cap=round,
  op/.style={font=\small,inner sep=1.6pt},
  lab/.style={font=\scriptsize,inner sep=0pt}
]
\node[op] (v1) at (0,0)        {$\alpha(v_1)$};
\node[op] (v2) at (-1.45,-.95) {$\alpha(v_2)$};
\node[op] (v3) at (1.45,-.95)  {$\alpha(v_3)$};
\node[op] (v4) at (2.90,-1.90) {$\alpha(v_4)$};

\draw (v1) -- node[lab,right=2pt,pos=.65] {$p$} (0,.60);

\foreach \a/\b/\ang in {v1/v2/-90,v1/v3/90,v3/v4/90}{
  \draw (\a) -- (\b);
  \coordinate (m) at ($(\a.center)!.5!(\b.center)$);
  \node[lab] at ($(m)!7pt!\ang:(\b.center)$) {$h$};
}

\foreach \a/\dx in {v2/-.85,v2/0,v2/.85,v3/-.85,
                    v4/-.85,v4/0,v4/.85}{
  \coordinate (leaf) at ($(\a.center)+(\dx,-.81)$);
  \draw (\a) -- (leaf);
  \coordinate (m) at ($(\a.center)!.66!(leaf)$);
  \node[lab] at ($(m)!5.5pt!90:(leaf)$) {$i$};
}
\end{tikzpicture}
\]
Such a transferred structure can depend on the chosen retract. Formality asks whether there is a choice of minimal model in which all higher operations vanish. The complex controlling this question is the convolution dg Lie algebra twisted by $\alpha_*$,
\[
\mathfrak g_{H}^{\alpha_*} \coloneqq \left(\mathfrak g_{H},\partial_{\alpha_*},[-,-] \right) \ ,
\qquad
\partial_{\alpha_*} \coloneqq  \partial+[\alpha_*,-] \ .
\]
Here we identify $\alpha_*$ with its degree-one, weight-one twisting cochain, given by $sv\mapsto\alpha_*(v)$ and zero in higher weights. The obstruction sequences theory of \cite{CE24b} then leads to the following result.

\begin{proposition}[Non-formality criterion]
\label{prop:successive-formality-obstruction}
Let $r\geqslant2$. If $\beta^{(k)}=0$ for $2\leqslant k<r$, then $\partial_{\alpha_*}\left(\beta^{(r)}\right) = 0$, and if the class
\[
\vartheta_r
\coloneqq \left[\beta^{(r)} \right]\in H^1\left(\mathfrak g_{H}^{\alpha_*}\right)
\]
is nonzero, then $(A,\alpha)$ is not formal.
\end{proposition}

\begin{proof}
This follows from \cite[Proposition 2.18 and 2.33]{CE24b} in the colored setting.
\end{proof}

\noindent For $r=2$, the class $\vartheta_2$ coincides precisely with Dimitrova's canonical class \cite[Theorem~4.2]{Dim12}, in the colored formulation of \cite[Section~2.7]{BSV24}.

\section{\textcolor{bordeau}{Locally constant prefactorization algebras}}\label{sec:locally-constant}

\noindent We now apply the obstruction theory of the previous section to locally constant prefactorization algebras on $\R^m$. We recall the two-dimensional Poisson bracket of \cite{BSV24} and extend its construction to higher dimensions using linking spheres. Under suitable vanishing assumptions, this yields a Poisson bracket of degree $1-m$ on cohomology whose non-vanishing implies non-formality.

\subsection{Prefactorization operad and algebras}
We start by recalling the definition of a prefactorization algebra on a smooth manifold, referring to \cite{CG17,CG21} for further details. 

\begin{definition}[Prefactorization operad] 
    Let $M$ be a smooth manifold of dimension $m$. Let $\mathcal D_M$ denote the set of open subsets of $M$ diffeomorphic to $\mathbb R^m$. The \emph{prefactorization operad} $\Pfac_M$ is the $\mathcal D_M$-colored operad defined as follows. For $\bbD=(D_1,\ldots,D_n)$ and $D\in\mathcal D_M$, set 
    \[ 
        \Pfac_M\! (D, \bbD) = \begin{cases} \mathbb{R} \bigl[\op{D}{\bbD}\bigr], & D_i\subseteq D \text{ for all }i \text{ and } D_i\cap D_j=\varnothing \text{ for }i\neq j,\\ 0,&\text{otherwise}, \end{cases} 
    \] 
    where $\mathbb{R}[\op{D}{\bbD}]$ is the cochain complex concentrated in degree zero generated by $\op{D}{\bbD}$. The symmetric group action is given by 
    \[ 
        \op{D}{\bbD_\sigma} = \op{D}{(D_{\sigma(1)},\ldots,D_{\sigma(n)})} \ , 
    \] 
    and the operadic composition is defined by 
    \[ 
        \gamma\left( \op{D}{\bbD}; \op{D_1}{\bbD_1},\ldots,\op{D_n}{\bbD_n} \right) = \op{D}{(\bbD_1,\ldots,\bbD_n)} \ .
    \]  
    The unit $\mathbf{1} : I \to \Pfac_M $  and the augmentation $\epsilon:\Pfac_M\to I$ are respectively given by 
    \[ 
        1_D=\op{D}{D}\in\Pfac_M (D,D)  \quad \mbox{and} \quad  \epsilon\left(\op{D}{D'}\right) = \begin{cases} 1_D,&D'=D,\\ 0,&\text{otherwise}. \end{cases}
    \] 
\end{definition}

\begin{definition}[Prefactorization algebra] Let $\Pfac_M$ be the
prefactorization operad over a smooth manifold $M$. A \emph{prefactorization algebra} on $M$ is an algebra over $\Pfac_M$. The prefactorization algebra $\mathfrak F$ is called \emph{locally constant} if, for every inclusion of disks $D\subseteq D'$, the structure map
\[
    \mathfrak F\!\left(\op{D'}{D}\right)\colon
    \mathfrak F(D)\xlongrightarrow{\simeq}\mathfrak F(D')
\]
is a quasi-isomorphism.
\end{definition}

\subsection{Locally constant prefactorization algebras on $\mathbb{R}^m$}
\label{subsec:locally_constant}
Let $\mathfrak F$ be a locally constant prefactorization algebra on $\R^m$,  $m\geqslant1$.  Let us fix a strong deformation retract $(i,p,h)$ of $\mathfrak F$ onto $\HF$. Let $(H, \beta)$ be a transferred structure on the latter given by Theorem \ref{thm:homotopy-transfer}. By definition, the disk inclusions induce isomorphisms on cohomology:
\[
    s_D\coloneqq\HF\!\left(\op{\R^m}{D}\right)
    \colon\HF(D)\overset{\cong}{\longrightarrow}H \ ,
\]
identifying the cohomology at every disk with $H\coloneqq\HF(\R^m)$. As observed in \cite{BSV24}, the isomorphisms $s_D$ can be used to improve the general choice of a strong deformation retract, and obtain a transferred structure $\tb$ on the constant family $H$.
\begin{lemma}[{\cite[Lemma~3.6]{BSV24}}]
    The maps
    \[
        \widetilde i_D=i_Ds_D^{-1} \ ,\qquad
        \widetilde p_D=p_{\R^m}\mathfrak F\!\left(\op{\R^m}{D}\right),
        \qquad
        \widetilde h_D=h_D-\widetilde i_D\widetilde p_Dh_D
    \]
    define a strong deformation retract from $\mathfrak F$ to the constant prefactorization algebra $H$.
\end{lemma}

\noindent There are then two different transferred structures on the constant family $H$: $\tb$, and the retract $\widehat\beta$ given by simply transporting $\beta$ along $s_D$, with a common strict structure
\[
    \mu_{\bbD}\coloneqq s_D 
    \HF\!\left(\op{D}{\bbD}\right)
    \bigl(s_{D_1}^{-1}\otimes\cdots\otimes s_{D_n}^{-1}\bigr) \ ,
    \qquad \bbD=(D_1,\ldots,D_n) \ .
\]
This strict structure is characterized as follows:
\begin{proposition}[{\cite[Proposition~3.5]{BSV24}}] \leavevmode
    \begin{enumerate}
        \item For $m\geqslant2$, the operations $\mu_{\bbD}$ depend on
        the number $n$ of input disks, but not on their positions. Writing $\mu_n$ for each such operation, the family $(\mu_n)_{n\geqslant0}$ defines a unital graded-commutative associative algebra structure on $H$.
        \item For $m=1$, the operations $\mu_{\bbD}$ also depend on the order of the intervals and determine instead a unital associative algebra structure on $H$.
    \end{enumerate}
\end{proposition}

\noindent The transferred structure $\tb$ along improved retract turns out to have nicer properties in what concerns its higher-weight operations. The following proposition extends parts of \cite[Proposition~3.7]{BSV24} to higher weights.

\begin{proposition}
    \label{prop:higher-weight-comparison}
    For every $r\geqslant2$, the component $\tb^{(r)}$
    \begin{enumerate}
        \item[(i)] It is constant in the outgoing disk, i.e. for all disks $D\subseteq\R^m$, we have
        \[
            \tb^{(r)}\!\left(t\left(\op{D}{\bbD};T_1,\ldots,T_n\right)\right)
            =\tb^{(r)}\!\left(t\left(\op{\R^m}{\bbD};T_1,\ldots,T_n\right)\right).
        \]
        \item[(ii)] 
        It vanishes on trees with unary root vertex, i.e. for all disks $D\subset D'$
        \[
            \tb^{(r)}\!\left(t\left(\op{D'}{D};T_0\right)\right)=0 \ .
        \]
        \item[(iii)] If $k\geqslant2$ and $\beta^{(r)}=0$ for $2\leqslant r<k$, then $\tb^{(r)}=0$ in the same range and for
        \[
            T=t\left(\op{D}{\bbD};T_1,\ldots,T_n \right)
        \] 
        of weight $k$, we have
        \[ 
            \tb^{(k)}(T)
            = \widehat\beta^{(k)}\!\left(
            t\left(\op{\R^m}{\bbD};T_1,\ldots,T_n\right)\right)
            -\sum_{\mathrm{wt}(T_i)=k-1}
            \mu_{\bbD}\circ_i
            \widehat\beta^{(k)}\!\left(t\left(\op{\R^m}{D_i};T_i\right)\right).
        \]
        Here $T_i$ has output $D_i$ and the sum contains at most one term.
    \end{enumerate}
   
\end{proposition}

\begin{proof}
    Properties~(i) and~(ii) follow from the tree formulas and the identities
    \[
        \widetilde p_D\mathfrak F\!\left(\op{D}{\bbD}\right)
        =p_{\R^m}\mathfrak F\!\left(\op{\R^m}{\bbD}\right),
        \qquad
        \widetilde p_{D'}\mathfrak F\!\left(\op{D'}{D}\right)
        \widetilde h_D
        =\widetilde p_D\widetilde h_D=0 \ . 
    \]  
    The case $k = 2$ of property (iii) corresponds to \cite[(3.18)]{BSV24}, which generalizes to higher weights as follows. Expanding the expression of $\widetilde h_D$ along the internal edges of a transferred tree gives a signed sum indexed by collections of cuts. If a collection produces a root component with $a$ vertices and non-root components with $b_1,\ldots,b_q$ vertices, the corresponding term is a composition of components of $\widehat\beta$ of weights
    \[
        a \ ,\quad b_1+1,\ldots,b_q+1 \ .
    \]
    Since $\widehat\beta^{(r)}=0$ if and only if $\beta^{(r)}=0$, the hypothesis implies that contributions in weights $2,\ldots,k-1$ vanish. In weight $k$, only two types of terms can remain: the empty cut, and a single cut separating a subtree of weight $k-1$ from the root vertex. Their contributions are
    \[
        \widehat\beta^{(k)}
        \left(t\left(\op{\R^m}{\bbD};T_1,\ldots,T_n\right)\right)
        \quad 
        \mbox{and}
        \quad
        -\mu_{\bbD}\circ_i
        \widehat\beta^{(k)}
        \left(t\left(\op{\R^m}{D_i};T_i\right)\right),
    \]
    proving the desired formula.
\end{proof}

\noindent In many of our arguments, we will be interested in invariants that depend not on the specific form of the components $\tb^{(k)}$, but only on their cohomology class. In principle, this would require one to show the invariant in question vanishes on any cochain of the form $\partial_{\mu} \chi$. It turns out that it is sufficient to only consider $\chi$ satisfying a certain normalization condition. This generalizes the discussion of \cite[Remark 3.8]{BSV24}.

\begin{definition}
    \label{def:normalized-cochain}
    A homogeneous cochain in $\mathfrak g_{H}^\mu$, supported in positive weights, is called \emph{normalized} if its value is independent of the outgoing disk and vanishes on trees with unary root vertex.
\end{definition}

\begin{lemma}
    \label{lem:normalized-coboundary}
    The differential $\partial_\mu$ preserves normalized cochains. Moreover, let $\alpha,\alpha'$ be normalized degree-one cocycles supported on weight $k\geqslant2$. If
    \[
        [\alpha]=[\alpha']
        \quad\text{in}\quad
        H^1\bigl(\mathfrak g_H^\mu\bigr) \ ,
    \]
    then there exists a normalized degree-zero cochain $\chi$, supported on weight $k-1$, such that
    \[
        \alpha'-\alpha=\partial_\mu\chi \ .
    \]
\end{lemma}

\begin{proof}
Let $\chi$ be a normalized cochain of degree $q$.
For $D\subsetneq D'$, let $T_1$ have output $D$, let $T_2$
be the same tree with output $D'$, and set
$W=t(\iota_D^{D'};T_1)$. Let us recall that
\[
    \partial_\mu\chi
    =(-1)^{q+1}\chi\circ d_{B\mathcal P}
     +\mu\star\chi-(-1)^q\chi\star\mu \ .
\]
In $\mu\star\chi$, only the cut below the new root contributes:
since $\mu_1=\operatorname{id}_H$, it gives
$(-1)^q\chi(T_1)$. Every term of $\chi\star\mu$ vanishes,
because $\chi$ is evaluated on a unary-root tree.
Likewise, the only surviving bar contraction is the one
immediately below the new root, which gives $T_2$ with sign $+1$.
Hence
\begin{equation}
    \label{eq:unary-root-differential}
    (\partial_\mu\chi)(W)
    =(-1)^q\bigl(\chi(T_1)-\chi(T_2)\bigr) \ .
\end{equation}
The right-hand side vanishes by normalization of $\chi$.
Thus $\partial_\mu\chi$ vanishes on unary-root trees, 
vanishing in weight one follows because $\partial_\mu$ raises
weight. The identity above uses only vanishing on unary-root
trees. Applying it to $\partial_\mu\chi$ and using
$\partial_\mu^2\chi=0$ proves independence of the outgoing disk.
Therefore $\partial_\mu\chi$ is normalized.

\medskip
\noindent
Now choose a degree-zero cochain $\xi$ of weight $k-1$ such that
$\partial_\mu\xi=\alpha'-\alpha$. For a tree $T$ with output $D$,
set
\[
    \eta(T)=
    \begin{cases}
    \xi\bigl(t(\op{\R^m}{D};T)\bigr) \ ,
        & D\subsetneq\R^m \ ,\\
    0 \ ,  & D=\R^m \ ,
    \end{cases}
\]
again ordering the added root first. Then $\eta$ has degree $-1$
and weight $k-2$. For $k=2$, we temporarily use the full
convolution complex, including the coaugmentation terms,
as in \cite[(2.39)]{BSV24}. Its differential still squares
to zero, and $\partial_\mu\eta$ has positive weight. Set $\chi=\xi+\partial_\mu\eta$, so that
$\partial_\mu\chi=\alpha'-\alpha$.
For a unary-root tree $T$ with output $D\subsetneq\R^m$,
expanding the differential gives
\[
    (\partial_\mu\xi)\bigl(t(\op{\R^m}{D};T)\bigr)
    =\xi(T)+(\partial_\mu\eta)(T)
    =\chi(T) \ .
\]
The added root contributes $\xi(T)$, and the remaining terms
give $\partial_\mu\eta(T)$.
The left-hand side is zero since $\alpha'-\alpha$ is normalized. If $T=t\left(\op{\R^m}{E};S\right)$ has output $\R^m$, all terms evaluating
$\eta$ on a tree with that output vanish. Only the root action
remains:
\[
    (\partial_\mu\eta)(T)=-\eta(S)=-\xi(T) \ .
\]
Thus $\chi$ vanishes on all unary-root trees.
Finally, applying \eqref{eq:unary-root-differential} to $\chi$,
whose differential is normalized, proves independence of the
outgoing disk. Hence $\chi$ is normalized.
\end{proof}

\subsection{The BSV Poisson bracket in dimension two}

Let $\mathfrak F$ be a locally constant prefactorization algebra
on $\R^2$, and write $H=\HF(\R^2)$.
Consider the first obstruction to formality $\vartheta_2$, as in Proposition~\ref{prop:successive-formality-obstruction}.
As explained in \cite[Section~3.4]{BSV24}, this obstruction
is associated with a concrete binary invariant,
whose construction we now recall. Let us choose disks
\[
    D=\phi(B_2(0)) \ ,\qquad
    \widetilde D=\phi(B_1(0)) \ ,
\]
where $\phi:\overline{B_2(0)}\to\R^2$ is a smooth embedding and $B_r(0)$ is the open ball of radius $r$. Inside $D$, choose four disks $D^1_1,D^1_2,D^2_1,D^2_2$ satisfying
\[
    D^1_1\cup D^1_2\subseteq D^2_1\cap D^2_2 \ ,
    \qquad
    (D^2_1\cup D^2_2)\cap\widetilde D=\varnothing \ .
\]
Let us assume that $D^2_1\cup D^2_2$ forms an annulus around $\widetilde D$, with $D^1_1$ and $D^1_2$ lying in the two distinct components of $D^2_1\cap D^2_2$. 
\begin{figure}[htbp]
    \centering
    \includegraphics[width=0.3\textwidth]{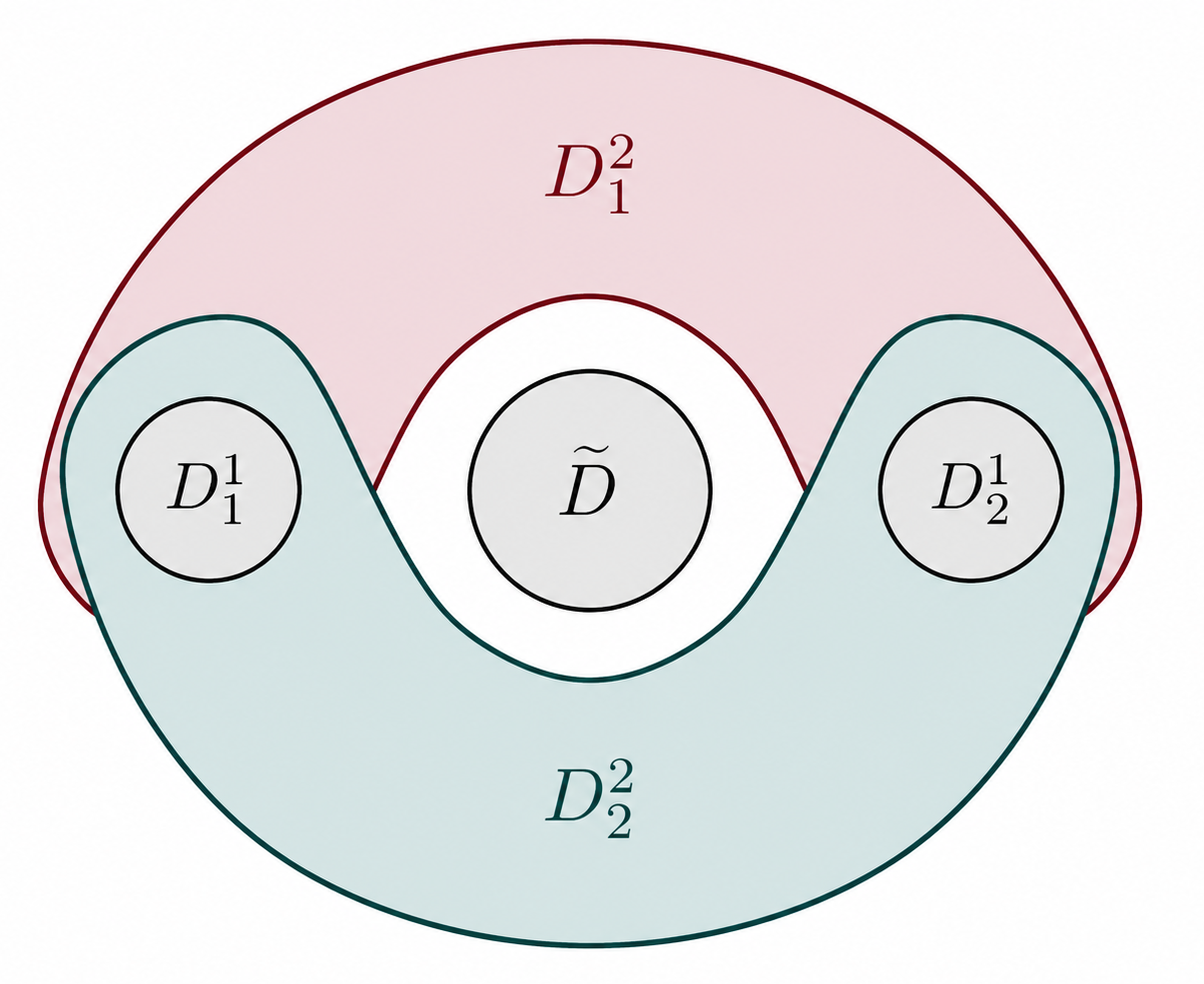}
    \caption{The disk configuration defining the two-dimensional bracket $L^{(2)}$.}
    \label{fig:cross-polytope-dim2}
\end{figure} 

\noindent Evaluating $\tb^{(2)}$ along its four inclusions gives the alternating sum
\[
    \begin{aligned}
        L^{(2)} \coloneqq {}&
        \tb^{(2)}
        \left(
        t\left(
        \iota^D_{(D^2_1,\widetilde D)},
        \iota^{D^2_1}_{D^1_1}
        \right)
        \right)
        -
        \tb^{(2)}
        \left(
        t\left(
        \iota^D_{(D^2_1,\widetilde D)},
        \iota^{D^2_1}_{D^1_2}
        \right)
        \right)
        \\
        &+
        \tb^{(2)}
        \left(
        t\left(
        \iota^D_{(D^2_2,\widetilde D)},
        \iota^{D^2_2}_{D^1_2}
        \right)
        \right)
        -
        \tb^{(2)}
        \left(
        t\left(
        \iota^D_{(D^2_2,\widetilde D)},
        \iota^{D^2_2}_{D^1_1}
        \right)
        \right).
    \end{aligned}
\]
Each summand is a weight-two operation, so $L^{(2)}:H^{\otimes2}\to H$ has degree $-1$. The alternating sum cancels the contributions of normalized coboundaries and $L^{(2)}$ depends only on $\vartheta_2$, by \cite[Proposition~3.10]{BSV24}. It is also independent of the disk configuration with the chosen orientation \cite[Lemma~3.11]{BSV24}. Together with the induced product on $H$, $L^{(2)}$ defines a $\mathbb P_2$-algebra structure \cite[Theorem~3.12]{BSV24}, and provides a concrete non-formality criterion
\[
    L^{(2)}\neq0
    \quad\Longrightarrow\quad
    \vartheta_2\neq0
    \quad\Longrightarrow\quad
    \mathfrak F\text{ is not formal}.
\]

\noindent Let us prove that the converse of the implication above does not hold with a counterexample. Consider the following dg algebra :
\[ 
    A = \Lambda(x,y,z) \ , \ |x| = |y| = |z| = 1 \ , \ dx = dy = 0 \ , \ dz = xy \ .
\]
This dg algebra admits a homotopy retract where $i\colon H \hookrightarrow A$ is the inclusion as the subspace
\[ 
    \langle1,x,y,xz,yz,xyz\rangle \ ,
\] 
the projection $p:A\to H$ sends $xy$ and $z$ to zero and the other basis elements to their classes, and we have $h(xy)=-z$ and $h=0$ on the other basis elements.  The dg algebra $A$ is a classic example of an algebra with a non-zero Massey product: products of degree one classes vanish, that is, $H^1 H^1 = 0$, so $[xy]=0$ and $[yy] = 0$, but there is a non-trivial triple product
\[ 
    [x,y,y] = [zy] = -[yz]
\]
that prevents $A$ from being formal. We define $\mathfrak{F}$ to be the \emph{constant} prefactorization algebra on $\R^m$, for $m \geqslant 1$, which assigns $\mathfrak F(D)=A$ for every disk $D$.
\begin{proposition}\label{prop:massey-counterexample}
    For every $m\geqslant1$, we have $\vartheta_2\neq0$. In dimension $m=2$, where the invariant $L^{(2)}$ is defined, one nevertheless has $L^{(2)}=0$.
\end{proposition}

\begin{proof}
    In dimension two, each of the four terms defining $L^{(2)}$ vanishes: the unary structure maps are identities, so the corresponding transferred operation contains the factor $hi=0$. We now prove that $\vartheta_2 \neq0$ for every $m\geqslant1$.
    Inside a fixed disk $D$, choose three disjoint disks $E_1,E_2,E_3$, a disk $D_{12}$ containing $E_1,E_2$ and disjoint from $E_3$, and a disk $D_{23}$ containing $E_2,E_3$ and disjoint from $E_1$.  In dimension one, choose the intervals in the order $E_1<E_2<E_3$. With root-first vertex order, set
    \[
        T_L=t\bigl(\op{D}{(D_{12},E_3)},\op{D_{12}}{(E_1,E_2)}\bigr) \ , \qquad
        T_R =t\bigl(\op{D}{(E_1,D_{23})},\op{D_{23}}{(E_2,E_3)}\bigr) \ .
    \]
    Evaluating $\tb^{(2)}$ on these trees and applying the resulting map to $[x]\otimes[y]\otimes[y]$, we get
    \[
        \tb^{(2)}(T_L)([x],[y],[y])=p\bigl(h(xy)y\bigr)=p(-zy)=[yz]\neq0 \ ,
    \]
    whereas
    \[
        \tb^{(2)}(T_R)([x],[y],[y])=0 \ ,
    \]
    so $\tb^{(2)}$ cannot be exact, as this would imply equality of these two terms.
\end{proof}

\subsection{Higher-dimensional Poisson bracket}
\label{subsec:higher-dimensional-poisson}
In this section, we extend the construction of $L^{(2)}$ to every dimension $m\geqslant2$. Let $\mathfrak F$ be a locally constant prefactorization algebra on $\R^m$, and write $H\coloneqq\HF(\R^m)$. Let us suppose that
\[
    \tb^{(2)}=\cdots=\tb^{(m-1)}=0 \ .
\]
Let us fix disks $D=\phi(B_2(0))$ and $\widetilde D=\phi(B_1(0))$ in $\mathbb{R}^m$ and $\mathscr D(D;\widetilde D)$ the poset of disks in 
\[
    D\setminus\overline{\widetilde D} \ ,
\] 
ordered by inclusion. These disks form a basis of contractible open sets. The classical covering lemma, see for instance \cite[Appendix~A.7, Lemma~A.1]{DKS23}, gives
\[
    H_{m-1}\bigl(N_\bullet\mathscr D(D;\widetilde D);\R\bigr)
    \cong H_{m-1}\bigl(D\setminus\overline{\widetilde D};\R\bigr)
    \cong\R \ ,
\]
since the complement is diffeomorphic to $S^{m-1}\times(1,2)$. We denote by $[S^{m-1}]$ the positive generator,
represented by an outward-oriented sphere surrounding $\widetilde D$.
\begin{definition}[Cross-polytope disk configurations]
    \label{def:cross-polytope-configuration}
    A \emph{cross-polytope disk configuration} $\mathcal C$ consists of disks $D^r_1,D^r_2\in\mathscr D(D;\widetilde D)$, for $1\leqslant r\leqslant m$, such that
    \[
        D^1_1\cap D^1_2=\varnothing \ ,
        \qquad
        D^{r-1}_1\cup D^{r-1}_2\subsetneq D^r_j
        \quad(2\leqslant r\leqslant m,\ j\in\{1,2\}) \ .
    \]
    Its associated cycle is
    \[
        z_{\mathcal C}\coloneqq
        \sum_{\varepsilon\in\{1,2\}^m}(-1)^{|\varepsilon|}
        \bigl[D^1_{\varepsilon_1}\subsetneq\cdots
        \subsetneq D^m_{\varepsilon_m}\bigr] \ ,
        \qquad
        |\varepsilon|\coloneqq\#\{r\mid\varepsilon_r=2\} \ .
    \]
    We call $\mathcal C$ \emph{positively oriented} if $[z_{\mathcal C}]=[S^{m-1}]$.
\end{definition}

\noindent For example, the following figure illustrates such a
configuration in dimension three. Choosing one disk at each
of the three levels gives the eight terms of $z_{\mathcal C}$.

\begin{figure}[htbp]
    \centering
    \includegraphics[width=.4\linewidth,trim=0 27bp 0 0,clip]{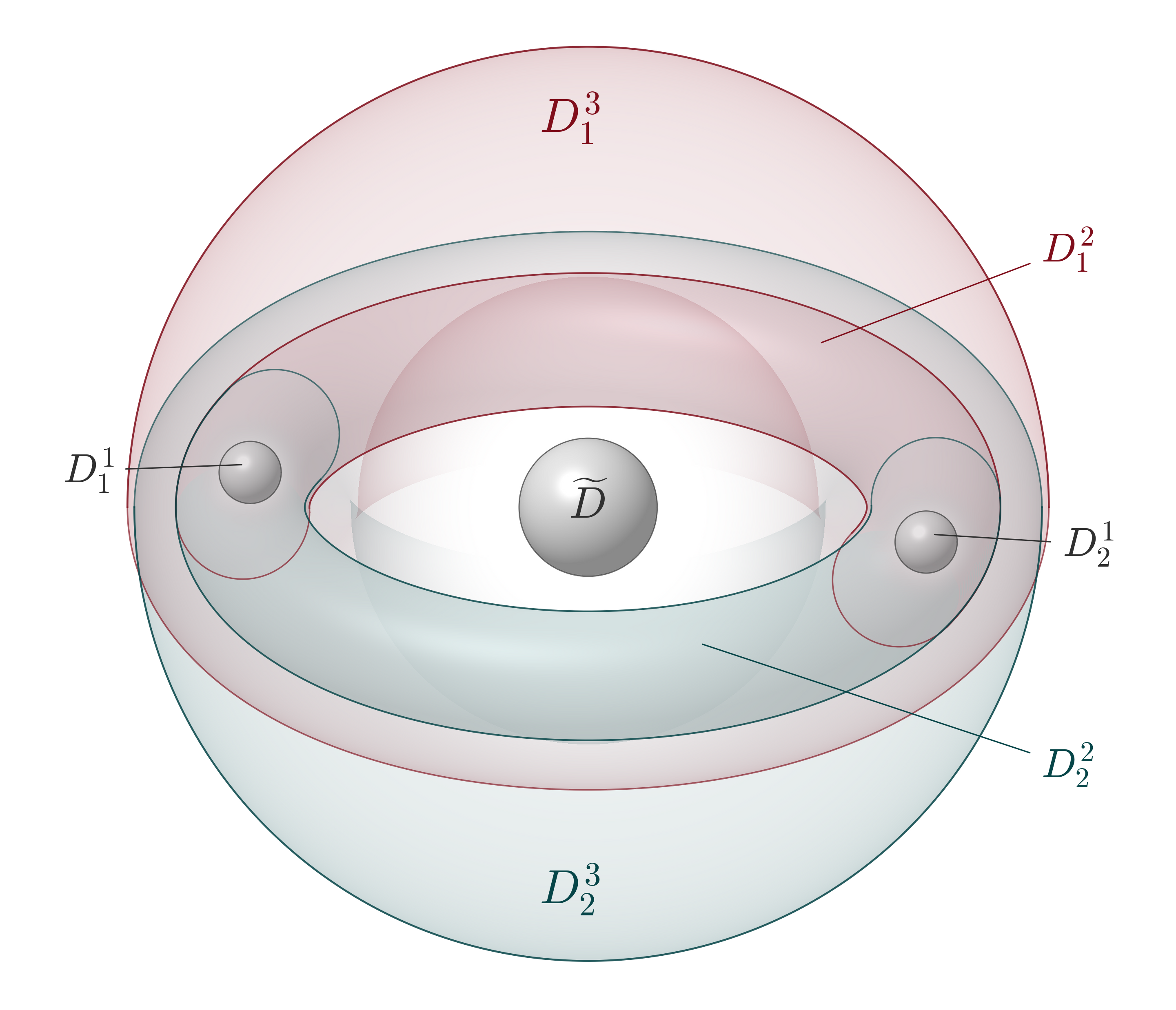}
    \caption{A cross-polytope disk configuration in dimension three.}
    \label{fig:cross-polytope-dim3}
\end{figure}

\noindent \noindent Each disk chain $D_1\subsetneq\cdots\subsetneq D_r$ in $\mathscr D(D;\widetilde D)$ determines a tree whose unary branch records these inclusions and whose binary root combines $D_r$ with $\widetilde D$ inside $D$:
\[
    T(D_1,\ldots,D_r)\coloneqq
    t\left(
    \op{D}{(D_r,\widetilde D)},
    \op{D_r}{D_{r-1}},\ldots,\op{D_2}{D_1}
    \right).
\]
We order the vertices from the root down the unary branch.

\begin{definition}[Higher-dimensional operation]
    \label{def:higher-poisson-invariant}
    For a positively oriented configuration $\mathcal C$,
    define the binary operation
    \[
        L^{(m)}_{\mathcal C}
        \coloneqq
        \sum_{\varepsilon\in\{1,2\}^m}(-1)^{|\varepsilon|}
        \tb^{(m)}\!\left(
            T(D^1_{\varepsilon_1},\ldots,D^m_{\varepsilon_m})
        \right).
    \]
    The sum runs over the $2^m$ choices of one disk at each
    level. Each summand is a binary operation on $H$
    of degree $1-m$.
\end{definition}

\noindent To compare configurations, we interpret the sum as the evaluation of a simplicial cochain on $z_{\mathcal C}$. Let $f\in\mathfrak g_H $ have degree $q$ and weight $r\geqslant1$. To define $\lambda(f)$ on a chain of disks, we evaluate $f$ on the corresponding tree as follows
\[
    \left\langle\lambda(f),
    [D_1\subsetneq\cdots\subsetneq D_r]\right\rangle
    \coloneqq f\left(T(D_1,\ldots,D_r)\right),
\]
for all disk chains in $\mathscr D(D;\widetilde D)$. Since a chain of $r$ disks is an $(r-1)$-simplex, $\lambda(f)$ is a simplicial $(r-1)$-cochain, with values in binary operations on $H$ of degree $q-r$. By linearity,  $f = \tb^{(m)}$ gives the alternating sum in Definition~\ref{def:higher-poisson-invariant}:
\[
    L^{(m)}_{\mathcal C} = 
    \left\langle\lambda \left(\tb^{(m)} \right ),z_{\mathcal C}\right\rangle.
\]

\begin{proposition}
    \label{prop:disk-chain}
    Let $f\in\mathfrak g_H^\mu$ have degree $q$, be supported on weight $r\geqslant1$, and vanish on trees with unary root. Then
    \begin{equation}
        \label{eq:disk-chain-differential}
        \lambda(\partial_\mu f)
        =(-1)^{q+r+1}\delta\lambda(f)
    \end{equation}
    and $\lambda(\tb^{(m)})$ is a simplicial cocycle. For a normalized degree-zero cochain $\chi$ of weight $m-1$, 
    \[
        \lambda\left(\tb^{(m)}+\partial_\mu\chi \right)
        -\lambda\left(\tb^{(m)}\right)=(-1)^m\delta\lambda(\chi) \ .
    \]
\end{proposition}

\begin{proof}
Evaluate $\partial_\mu f$ on the tree associated with $U_0\subsetneq\cdots\subsetneq U_r$, with internal edges numbered $k=0,\ldots,r-1$ from the root. Contracting edge $k$ deletes $U_{r-k}$, and its bar sign $(-1)^k$ gives the coefficient $(-1)^{q+1+k}$ in the convolution differential. The term $\mu\star f$ vanishes because it evaluates $f$ on a tree with unary root. The only term in $-(-1)^q f\star\mu$ deletes $U_0$, using $\mu_1=\operatorname{id}_H$. Moving $\mu$ past the $r$ odd vertex decorations contributes $(-1)^r$, giving the coefficient $(-1)^{q+r+1}$. Factoring out this coefficient leaves precisely the simplicial boundary signs, proving \eqref{eq:disk-chain-differential}. The two consequences follow by taking $(q,r)=(1,m)$ and $(q,r)=(0,m-1)$, respectively.
\end{proof}

\begin{proposition}
    \label{prop:higher-poisson-invariant}
    The operation $L^{(m)}_{\mathcal C}$ is independent of the positively oriented configuration $\mathcal C$, including the disks $D$ and $\widetilde D$. The resulting operation $L^{(m)}$ depends only on
    \[
        \vartheta_m\coloneqq\left[\tb^{(m)}\right]\in H^1(\mathfrak g_H^\mu) \ .
    \]
    In particular, $L^{(m)}\neq0$ implies $\vartheta_m\neq0$ and hence non-formality of $\mathfrak F$.
\end{proposition}

\begin{proof}
    For fixed $D,\widetilde D$, the cochain is closed by Proposition~\ref{prop:disk-chain}, so its pairing depends only on the positive linking class $[z_{\mathcal C}]$. To vary the disks, we also allow any outer disk $D$ containing the smoothly embedded closed disk $\overline{\widetilde D}$. Independence under variation of $D$ follows from the fact that $\tb^{(m)}$ is independent of the outgoing disk. For the inner disk, let $\widetilde D\subsetneq\widetilde D'$ and restrict the two corresponding cochains, denoted by $\lambda$ and $\lambda'$, to the disk poset of $D\setminus\overline{\widetilde D'}$. Define a simplicial $(m-2)$-cochain $\kappa$ by setting its value on a  $U_0\subsetneq\cdots\subsetneq U_{m-2}$ to be the evaluation of $\tb^{(m)}$ on the following tree:
    \begin{center}
        \begin{tikzpicture}[every node/.style={inner sep=1.5pt,font=\small}]
            \node (r) at (0,0) {$\op{D}{(U_{m-2},\widetilde D')}$};
            \draw (r) -- ++(0,.75);
            \node (u1) at (-1.2,-1.3) {$\op{U^{m-2}}{U^{m-3}}$};
            \node (u2) at (-1.2,-2.6) {$\op{U^{1}}{U^{0}}$};
            \node (e) at (1.3,-1.3) {$\op{\widetilde D'}{\widetilde D}$};
            \draw (r) -- (u1);
            \draw (r) -- (e);
            \draw[densely dotted] (u1) -- (u2);
            \draw (u2) -- ++(0,-.8);
            \draw (e) -- ++(0,-.9);
        \end{tikzpicture}
    \end{center}
    We set the sign by ordering the root first, then the vertices on the first branch, and finally the vertex on the second branch. Evaluating $\partial_\mu\tb^{(m)}=0$ on the analogous tree with one more unary vertex on the left branch tree gives
    \[
        0=(-1)^{m-1}\bigl(\delta\kappa+\lambda-\lambda'\bigr) \ .
    \]
    The first-branch terms give the coboundary, while those on the second branch compare the two inner disks. Consequently, $\lambda'-\lambda=\delta\kappa$, proving independence on the choice of $\widetilde D$.
    \medskip
    
    \noindent For a pair arbitrary configurations, choose concentric round disks $E\subset B$ with $E$ containing the closures of both $\widetilde D$ disks. We can enlarge the outgoing disks to $B$ and compare each of the other disks with $E$, using a positive cycle in $B\setminus\overline E$, proving independence of all disk choices. Finally, two normalized weight-$m$ cocycles representing the same class differ by $\partial_\mu\chi$ with $\chi$ normalized of degree zero and weight $m-1$, by Lemma~\ref{lem:normalized-coboundary}. By \eqref{eq:disk-chain-differential}, their simplicial cochains differ by a coboundary and have the same pairing with the linking class. In particular, $\vartheta_m=0$ implies $L^{(m)}=0$, so the non-formality assertion follows from Proposition~\ref{prop:successive-formality-obstruction}.
\end{proof}

\begin{theorem}
    \label{thm:higher-poisson-bracket}
    The operation $L^{(m)}$ defines a Poisson bracket of degree $1-m$ on the associative, unital, and commutative algebra $(H,\mu_2)$. Explicitly, it satisfies
    \begin{enumerate}[label=(\roman*)]
        \item
    $L^{(m)}\tau=(-1)^mL^{(m)}$ \ ,

    \smallskip
        \item
        $ L^{(m)}(\operatorname{id}\otimes\mu_2)
        = \mu_2(L^{(m)}\otimes\operatorname{id})   +\mu_2(\operatorname{id}\otimes L^{(m)})
        (\tau\otimes\operatorname{id}) \ ,$
        \smallskip
        \item
        $L^{(m)}(\operatorname{id}\otimes L^{(m)})+L^{(m)}(\operatorname{id}\otimes L^{(m)})\tau_{(123)} +L^{(m)}(\operatorname{id}\otimes L^{(m)})\tau_{(132)}
        =0 \ . $
    \end{enumerate}
    Here all permutations act through the symmetric braiding. 
\end{theorem}

\begin{proof}
    For $n=2,3$, consider product cells
    \[
        \sigma=\prod_{i=1}^n
        [U_i^0\subsetneq\cdots\subsetneq U_i^{p_i}] \ ,
    \]
    with pairwise disjoint largest disks. They form a multisimplicial set whose diagonal is the nerve of the poset of ordered $n$-tuples of pairwise disjoint disks. By the Eilenberg--Zilber theorem, their total chain complex, with boundary operator
    \[ 
        \partial\sigma=\sum_i(-1)^{p_1+\cdots+p_{i-1}}\partial_i\sigma \ ,
    \]
    where $\partial_i$ denotes the simplicial boundary of the $i$-th factor, computes the homology of this nerve. The same covering argument as we used for disk chains, applied to the product opens $U_1\times\cdots\times U_n$, identifies the homology of the nerve with the homology of the configuration space $H_*(\operatorname{Conf}_n(\R^m);\R)$. For a homogeneous cochain $f\in\mathfrak g_H$ of weight $r\geqslant1$, define a cochain on this complex of product cells by setting
    \[
        \langle\lambda_n(f),\sigma\rangle
        =(-1)^{\sum_{i<j}p_ip_j}f(T_\sigma) \ ,
        \qquad \sum_i p_i=r-1 \ .
    \]
    Here $T_\sigma$ is the following rooted tree:
    \begin{center}
        \begin{tikzpicture}[every node/.style={inner sep=1.5pt,font=\small}]
            \node (r) at (0,0) {$\op{\R^m}{(U_1^{p_1},\ldots,U_n^{p_n})}$};
            \draw (r) -- ++(0,.75);
            \node (a1) at (-2.3,-1.3) {$\op{U_1^{p_1}}{U_1^{p_1-1}}$};
            \node (a2) at (-2.3,-2.8) {$\op{U_1^{1}}{U_1^{0}}$};
            \node (b1) at (2.3,-1.3) {$\op{U_n^{p_n}}{U_n^{p_n-1}}$};
            \node (b2) at (2.3,-2.8) {$\op{U_n^{1}}{U_n^{0}}$};
            \node at (0,-1.3) {$\cdots$};
            \draw (r) -- (a1);
            \draw (r) -- (b1);
            \draw[densely dotted] (a1) -- (a2);
            \draw[densely dotted] (b1) -- (b2);
            \draw (a2) -- ++(0,-.8);
            \draw (b2) -- ++(0,-.8);
        \end{tikzpicture}
    \end{center}
    We fix the sign by ordering the root first, then the branches in input order, each from root to leaf. Note that reordering two branches in the bar construction costs the sign $(-1)^{p_ip_j}$ of exchanging the corresponding factors of a product cell, so these cochains are equivariant under permutations of the inputs. 
    \medskip
    
    \noindent Now let $f$ be a cochain of degree $q$, weight $r$, and vanish on trees with unary root, and evaluate $\partial_\mu f$ on $T_\sigma$, where $\sum_ip_i=r$. As in Proposition~\ref{prop:disk-chain}, the term $\mu\star f$ vanishes, and there is a term  deleting each $U_i^a$, come from contractions or from $f\star\mu$. All of these terms carry the sign $(-1)^{q+1+p_1+\cdots+p_i+a}$, and accounting for the fact that $\lambda_n$ changes by $(-1)^{r-p_i}$ when $p_i$ decreases by one, the final sign is$(-1)^{q+r+1}$, so we get:
    \begin{equation}
        \label{eq:compatibility-lambda-differentials}
        \lambda_n(\partial_\mu f)=(-1)^{q+r+1}\delta\lambda_n(f) \ .
    \end{equation}
    In particular, at weight $r=m$, $\lambda_n(\tb^{(m)})$ is closed. To prove the symmetry statement \emph{(i)}, we look at arity $n=2$, where pairing with a positive linking class $[z_{\mathcal C}]$ gives the operation $L^{(m)}$. Exchanging the inputs induces the antipodal map on $\operatorname{Conf}_2(\R^m)\simeq S^{m-1}$, which has degree $(-1)^m$, so the symmetry of $L^{(m)}$ (with that sign) follows from equivariance.
    \medskip
    
    \noindent For the Leibniz identity, fix disjoint round disks $E_2,E_3$ for inputs $2,3$, and let $c_{1\mid23},c_{12},c_{13}$ be products of a positively oriented cross-polytope cycle for input $1$ with $[E_2]\times[E_3]$, in which input $1$ surrounds both fixed disks, only $E_2$, or only $E_3$, respectively. They satisfy
    \[
        [c_{1\mid23}]=[c_{12}]+[c_{13}]
        \quad\text{in }H_{m-1}(\operatorname{Conf}_3(\R^m);\R) \ .
    \]
    Choose these cycles together with intermediate disks: for $c_{1\mid23}$, a round disk $E$ containing $\overline{E_2}\cup\overline{E_3}$ and disjoint from all the `moving disks', that is, the disks appearing in the disk chain attached to input $1$. For $c_{12}$, we choose a disk $E'$ containing $\overline{E_2}$ and the moving disks, disjoint from $E_3$, and for $c_{13}$ we do the same, exchanging labels $2,3$. Evaluate $\partial_\mu\tb^{(m)}=0$ on the following trees, where $U^0\subsetneq\cdots\subsetneq U^{m-1}$ runs over the chains of the cycle, and sum.
    \begin{center}
        \begin{tikzpicture}[every node/.style={inner sep=1.5pt,font=\small}]
            \begin{scope}
                \node (r) at (0,0) {$\op{\R^m}{(U^{m-1},E)}$};
                \draw (r) -- ++(0,.75);
                \node (u1) at (-1.2,-1.3) {$\op{U^{m-1}}{U^{m-2}}$};
                \node (u2) at (-1.2,-2.6) {$\op{U^{1}}{U^{0}}$};
                \node (e) at (1.3,-1.3) {$\op{E}{(E_2,E_3)}$};
                \draw (r) -- (u1);
                \draw (r) -- (e);
                \draw[densely dotted] (u1) -- (u2);
                \draw (u2) -- ++(0,-.8);
                \draw (e) -- ++(-.5,-.9);
                \draw (e) -- ++(.5,-.9);
                \node at (0,-4.1) {for $c_{1\mid23}$};
            \end{scope}
            \begin{scope}[xshift=6.2cm]
                \node (r) at (0,0) {$\op{\R^m}{(E',E_3)}$};
                \draw (r) -- ++(0,.75);
                \node (e) at (-1,-1.1) {$\op{E'}{(U^{m-1},E_2)}$};
                \draw (r) -- (e);
                \draw (r) -- ++(.9,-1);
                \node (u1) at (-2,-2.2) {$\op{U^{m-1}}{U^{m-2}}$};
                \node (u2) at (-2,-3.3) {$\op{U^{1}}{U^{0}}$};
                \draw (e) -- (u1);
                \draw (e) -- ++(.6,-.9);
                \draw[densely dotted] (u1) -- (u2);
                \draw (u2) -- ++(0,-.7);
                \node at (0,-4.1) {for $c_{12}$};
            \end{scope}
        \end{tikzpicture}
    \end{center}
    In both cases, the terms deleting a disk of the unary chain are, as in the proof of Proposition~\ref{prop:higher-poisson-invariant}, the coboundary of a cochain evaluated on a cycle, so they vanish. For $c_{1\mid23}$, contracting the edge between the binary vertices contributes $(-1)^{m-1}\langle\lambda_3(\tb^{(m)}),c_{1\mid23}\rangle$, and detaching $\op{E}{(E_2,E_3)}$ contributes $(-1)^mL^{(m)}(\operatorname{id}\otimes\mu_2)$. For $c_{12}$, the contraction contributes $\langle\lambda_3(\tb^{(m)}),c_{12}\rangle$, and the term $\mu\star\tb^{(m)}$ contributes $-\mu_2(L^{(m)}\otimes\operatorname{id})$. The value on $c_{13}$ follows from that on $c_{12}$ by equivariance under the transposition of inputs $2$ and $3$. Hence
    \[
        \begin{aligned}
            \langle\lambda_3(\tb^{(m)}),c_{1\mid23}\rangle
            &=L^{(m)}(\operatorname{id}\otimes\mu_2) \ ,\\
            \langle\lambda_3(\tb^{(m)}),c_{12}\rangle
            &=\mu_2(L^{(m)}\otimes\operatorname{id}) \ ,\\
            \langle\lambda_3(\tb^{(m)}),c_{13}\rangle
            &=\mu_2(\operatorname{id}\otimes L^{(m)})
            (\tau\otimes\operatorname{id}) \ .
        \end{aligned}
    \]
    The homology relation between the three cycles now gives the Leibniz identity.
    \medskip 
    
    \noindent For Jacobi, the weight-$2m$ Maurer--Cartan equation reads
    \[
        \partial_\mu\tb^{(2m-1)}
        =-\tb^{(m)}\star\tb^{(m)} \ .
    \]
    Since the right-hand side vanishes on trees with unary root, the normalization calculation in the proof of Lemma~\ref{lem:normalized-coboundary}, applied in degree one, allows us to change $\tb^{(2m-1)}$ by a $\partial_\mu$-coboundary so that it also vanishes there. Denote the resulting cochain by $b$, which will satisfy the same equation above. Calculating the differential, on every $T_\sigma$, each cut has an entirely unary component, so $\lambda_3(b)$ is closed by \eqref{eq:compatibility-lambda-differentials}.
    \medskip
    
    \noindent Let $C_{23}=z_U\times z_V\times[E_3]$ be the product of two positively oriented cross-polytope cycles, $z_U$ around a round disk $E$ and $z_V$ around a round disk $E_3$ inside $E$, with the outer cycle first. Geometrically, it is represented by
    \[
        (u,v)\longmapsto(Ru,\rho v,0) \ ,
        \qquad (u,v)\in S^{m-1}\times S^{m-1} \ ,\quad R>2\rho>0 \ .
    \]
    Evaluate $\partial_\mu b+\tb^{(m)}\star\tb^{(m)}=0$ on the trees obtained by grafting the two binary-rooted chains, the inner one into the second input of the outer one, as follows:
    \begin{center}
        \begin{tikzpicture}[every node/.style={inner sep=1.5pt,font=\small}]
            \node (r) at (0,0) {$\op{\R^m}{(U^{m-1},E)}$};
            \draw (r) -- ++(0,.75);
            \node (u1) at (-1.4,-1.3) {$\op{U^{m-1}}{U^{m-2}}$};
            \node (u2) at (-1.4,-2.6) {$\op{U^{1}}{U^{0}}$};
            \node (e) at (1.6,-1.3) {$\op{E}{(V^{m-1},E_3)}$};
            \node (v1) at (0.7,-2.6) {$\op{V^{m-1}}{V^{m-2}}$};
            \node (v2) at (0.7,-3.9) {$\op{V^{1}}{V^{0}}$};
            \draw (r) -- (u1);
            \draw (r) -- (e);
            \draw[densely dotted] (u1) -- (u2);
            \draw (u2) -- ++(0,-.8);
            \draw (e) -- (v1);
            \draw (e) -- ++(.8,-1);
            \draw[densely dotted] (v1) -- (v2);
            \draw (v2) -- ++(0,-.8);
        \end{tikzpicture}
    \end{center}
    We fix the sign by ordering the outer root and its unary branch first, followed by the inner root and its unary branch. The terms deleting a disk of either unary chain are, as before, the coboundary of a cochain on products of chains, evaluated on the product cycle $z_U\times z_V$, so they vanish. The remaining terms contract or cut the edge between the binary vertices, giving
    \[
        \langle\lambda_3(b),C_{23}\rangle
        =(-1)^{m-1}L^{(m)}(\operatorname{id}\otimes L^{(m)}) \ .
    \]
    once we account for signs; the relative contraction sign is canceled by the factor $(-1)^{(m-1)^2}$ in $\lambda_3$, whereas the cut term which carries sign $(-1)^m$. 
    \medskip
    
    \noindent We then study the relation between this geometric cycle and its cyclic translates. In the unit sphere of the linear subspace
    \[
        \{(x_1,x_2,x_3)\in(\R^m)^3:x_1+x_2+x_3=0\} \ ,
    \]
    remove small disjoint tubular neighborhoods of the three collision spheres. The remaining compact oriented manifold has three boundary components. Around $x_2=x_3$, translation to $x_3=0$ identifies that boundary component with
    \[
        (u,v)\mapsto(Ru+\tfrac\rho2v,\rho v,0) \ ,
        \qquad R>2\rho>0 \ .
    \]
    Removing $\tfrac\rho2v$ gives a collision-free homotopy to $C_{23}$. Cyclic permutations preserve the ambient orientation and permute the boundary components, so the three comparison signs agree. Hence
    \[
        [C_{23}]+[(123)_*C_{23}]+[(132)_*C_{23}]=0 \ .
    \]
    Pairing with the closed cochain $\lambda_3(b)$ and using the preceding evaluation on $C_{23}$ proves Jacobi.  Lastly, to prove unitarity, the Leibniz identity applied to $1\cdot1$ gives $L^{(m)}(x,1)=0$, and symmetry gives the same in the first input.
\end{proof}

\section{\textcolor{bordeau}{Factorization envelopes}}
\label{subsec:factorization-envelopes}

\noindent Bruinsma, Schenkel and Vicedo studied the factorization envelope of a Lie algebra on $\R^m$ and computed its Poisson bracket in dimension two \cite[Section~4]{BSV24}. We extend this calculation to higher dimensions and show that, for $m\geqslant2$, the envelope is formal if and only if the Lie algebra is abelian. The same characterization holds for compactifications along a nonempty closed manifold. For envelopes of $L_\infty$-algebras, we generalize this result by proving that a nonzero induced Lie bracket on cohomology implies non-formality.

\subsection{Factorization envelopes of a Lie algebra}
Throughout this subsection $\mathfrak g$ is a finite-dimensional real or complex Lie algebra concentrated in degree zero, and we take differential forms over the same field. We apply the factorization envelope construction of \cite[Section~3.6]{CG17} to the local Lie algebra $\mathfrak g\otimes\Omega^\bullet_{\R^m}$. Its value on a disk $D$ is the Chevalley--Eilenberg chain complex
\[
    \mathfrak U^{(m)}_{\mathfrak g}(D)
    =
    \left(\Sym\bigl((\mathfrak g\otimes\Omega_c^\bullet(D))[1]\bigr),
          d_{\mathrm{dR}}[1]+d_{\mathrm{CE}}\right).
\]
The structure maps extend compactly supported forms by zero and multiply the resulting monomials. From now on, we leave these extension-by-zero maps implicit in our formulas.

\medskip

\noindent We use the retract of \cite[Section~4.1]{BSV24}. For every disk $D$, a compactly supported top form $\omega_D$ of integral one represents the generator of $H_c^m(D)$. There is a homotopy $k_D$ on compactly supported forms satisfying
\[
    d_{\mathrm{dR}}k_D(\omega)+k_D(d_{\mathrm{dR}}\omega)
       =\omega-\omega_D\int_D\omega \ ,
    \qquad k_D(\omega_D)=0 \ ,\qquad k_D^2=0 \ .
\]
On linear generators, the inclusion, projection and homotopy for the shifted de Rham complex are
\[
    i_D(a)=a\otimes\omega_D \ ,\qquad
    p_D(a\otimes\omega)=a\int_D\omega \ ,\qquad
    h_D(a\otimes\omega)=a\otimes k_D(\omega) \ .
\]
The maps $i_D$ and $p_D$ extend as a product, and for the homotopy $h_D$ we use the symmetrized tensor trick. 
We perturb the differential by $d_{\mathrm{CE}}$, giving perturbed maps
\[
    \widehat p_D=\sum_{j\geqslant0}p_D(d_{\mathrm{CE}}h_D)^j \ ,
    \qquad
    \widehat h_D=\sum_{j\geqslant0}h_D(d_{\mathrm{CE}}h_D)^j \ .
\]
Both sums are finite on every monomial, and moreover $d_{\mathrm{CE}}i_D=0$ because a wedge of two top forms vanishes. The inclusion $i_D$ and the zero differential on $H\coloneqq\Sym(\mathfrak g[1-m])$ therefore give an objectwise retract:
\[
    \begin{tikzcd}[column sep=large]
        H \arrow[r,shift left=.7ex,"i_D"] &
        \mathfrak U^{(m)}_{\mathfrak g}(D)
        \arrow[l,shift left=.7ex,"\widehat p_D"]
        \arrow[loop right,"\widehat h_D"]
    \end{tikzcd}
\]
As in Section~\ref{subsec:locally_constant}, $\widehat\beta$ denotes the transfer along these retracts, transported to the constant family $H$, and $\tb$ denotes the transfer along the improved retracts.

\begin{theorem}[{\cite[Section~4]{BSV24}}]
    \label{thm:factorization-envelope-low-dimensional} \leavevmode
   \begin{enumerate}
       \item For $ m = 1$, the strict transferred product on $\Sym(\mathfrak g)$ is the Gutt star product. All higher components vanish, and $\mathfrak U^{(1)}_{\mathfrak g}$ is formal.
       \item In every dimension $m\geqslant2$, the strict transferred product on $\Sym(\mathfrak g[1-m])$ is the usual graded commutative product. 
       \item For $m=2$, $L^{(2)}$ is the degree $-1$ Poisson bracket characterized on generators by 
       \[L^{(2)}(a,b)=[a,b] \ .\]
   \end{enumerate}
\end{theorem}

\noindent We now generalize these calculations to higher dimensions. Analyzing the symmetric degree limits the weights in which transferred operations can appear, giving the following results.

\begin{proposition}
    \label{prop:factorization-envelope-weight-vanishing}
    Let $m\geqslant2$. A transferred component $\widehat\beta^{(r)}$ can be non-zero only if
    \[
        r=1+q(m-1) \ ,
    \]
    for some $q\in\mathbb Z_{\geqslant0}$. The same statement holds for $\widetilde\beta^{(r)}$ and in particular we have
    \[
        \widehat\beta^{(2)} = \cdots = \widehat\beta^{(m-1)} = 0 \ ,
        \qquad
        \widetilde\beta^{(2)} = \cdots = \widetilde\beta^{(m-1)} = 0 \ .
    \]
\end{proposition}

\begin{proof}
A monomial of symmetric degree $n$ in $H$ has cohomological degree $n(m-1)$. For inputs containing $n$ generators in all, a weight-$r$ operation therefore has output degree $n(m-1)+1-r$. In $\widehat\beta^{(r)}$, each occurrence of $d_{\mathrm{CE}}$ replaces two generators by one. All the other maps preserve the number of generators. A nonzero term with $q$ occurrences must therefore satisfy
\[
    n(m-1)+1-r=(n-q)(m-1) \ .
\]
This gives $r=1+q(m-1)$, and the same degree-counting argument applies to $\tb^{(r)}$.
\end{proof}

\begin{proposition}\label{prop:abelian-formal}
    Let $m\geqslant2$. If $\mathfrak{g}$ is abelian, then $\widehat\beta^{(r)} = \tb^{(r)} = 0$ for every $r\geqslant2$, and thus $\mathfrak U\mathfrak g^{\R^m}$ is formal.
\end{proposition}
\begin{proof}
    As in the proof of the previous proposition, in any nonzero contribution to $\widehat{\beta}$ we must have $r-1 = q(m-1)$ where $q$ is the number of occurrences of $d_{\mathrm{CE}}$, so when $r \geqslant 2$ we must have $q \geqslant 1$, which implies vanishing since $d_{\mathrm{CE}}=0$ when $\mathfrak{g}$ is abelian. The analogous statement for $\tb$ follows then from \cref{prop:higher-weight-comparison}.
\end{proof}

\begin{corollary}
    \label{cor:factorization-envelope-weight-m}
    For $m\geqslant2$ and a tree $T$ with $m$ vertices, the only possibly nonzero summand of $\widehat\beta^{(m)}(T)$ has $h$ on each internal edge and $p\,d_{\mathrm{CE}}h$ on the outgoing edge.
\end{corollary}

\begin{proof}
    By Proposition~\ref{prop:factorization-envelope-weight-vanishing}, there is exactly one occurrence of $d_{\mathrm{CE}}$. If it lies on an internal edge $e$, let $T_e$ be the subtree below $e$, with $s$ vertices. Before reaching $d_{\mathrm{CE}}$, there are at most $s\leqslant m-1$ homotopies: one on $e$ and $s-1$ on the internal edges of $T_e$. The two forms paired by $d_{\mathrm{CE}}$ start in total degree $2m$, so their wedge can be nonzero only after at least $m$ homotopies, thus $d_{\mathrm{CE}}$ must lie on the outgoing edge.
\end{proof}

\noindent For the remainder of this section, let us assume $m\geqslant2$. We now extend the computation of \cite[Lemma~4.8 and Theorem~4.10]{BSV24} to weight $m$. Let
\[
\{-,-\}_{(1-m)}
\colon
H^{\otimes2}\longrightarrow H
\]
be the degree $1-m$ Poisson bracket determined by its values on linear monomials
\[
\{a, b\}_{(1-m)}
=
[a,b] \ ,
\qquad
a,b\in\mathfrak g \ .
\]
We consider a disk arrangement of the type used in the definition of $L^{(m)}$, given by a disk chain
\[
D^1\subsetneq D^2\subsetneq\cdots\subsetneq D^m
\]
with $D^m\subset D$ and $D^m\cap\widetilde D=\varnothing$. Denote by
\[
T(D^\bullet;\widetilde D)
\coloneqq
t\left(
\iota^D_{(D^m,\widetilde D)},
\iota^{D^m}_{D^{m-1}},
\ldots,
\iota^{D^2}_{D^1}
\right)
\]
the corresponding linear $m$-vertex tree, namely $T(D^1,\ldots,D^m)$ in the notation of Section~\ref{subsec:higher-dimensional-poisson}. Corollary~\ref{cor:factorization-envelope-weight-m} then gives
\begin{equation}
\label{eq:factorization-envelope-linear-weight-m}
\begin{aligned}
\widetilde\beta^{(m)}
\bigl(T(D^\bullet;\widetilde D)\bigr)
={}&
p_{\R^m}d_{\mathrm{CE}}h_{\R^m}
\circ
\mathfrak U\mathfrak g^{\R^m}
\left(
\iota^{\R^m}_{(D^m,\widetilde D)}
\right)
\circ
\left(
h_{D^m}\cdots h_{D^2}i_{D^1}
\otimes i_{\widetilde D}
\right)
\\
&-
\mu_2\circ_1
\left(
p_{\R^m}d_{\mathrm{CE}}h_{\R^m}
h_{D^m}\cdots h_{D^2}i_{D^1}
\right).
\end{aligned}
\end{equation}
We now set
\[
\alpha_{D^\bullet}
\coloneqq
k_{D^m}
k_{D^{m-1}}
\cdots
k_{D^2}
(\omega_{D^1})
\in
\Omega_c^1(D^m) \ .\]

\begin{lemma}
\label{lem:factorization-envelope-linear-weight-m}
For all $A,B\in H$, the operation on the linear tree $T(D^\bullet;\widetilde D)$ is given by
\[
\widetilde\beta^{(m)}
\bigl(T(D^\bullet;\widetilde D)\bigr)
(A\otimes B)
=
\varphi^{(m)}_0
\bigl(T(D^\bullet;\widetilde D)\bigr)
\{A,B\}_{(1-m)} \ ,
\]
that is, by the Poisson structure multiplied by a numerical factor
\[
\begin{aligned}
\varphi^{(m)}_0
\bigl(T(D^\bullet;\widetilde D)\bigr)
\coloneqq
\frac12\int_{\R^m}
\Big(
&
k_{\R^m}(\alpha_{D^\bullet})
\wedge
(\omega_{\widetilde D}+\omega_{\R^m})
 +
(-1)^m
k_{\R^m}(\omega_{\widetilde D})
\wedge
\alpha_{D^\bullet}
\Big)
\end{aligned}
\]
that is independent of $\mathfrak{g}$. 
\end{lemma}

\noindent For $m=2$, this is precisely \cite[(4.44)]{BSV24}, where the numerical factor function is denoted simply $\varphi_0$.

\begin{proof}
    Let us first consider the case where $A,B$ are linear, represented by $a,b\in\mathfrak g$, respectively. The $m-1$ internal homotopies transform
    \[
    a\otimes\omega_{D^1}
    \quad\text{into}\quad
    a\otimes\alpha_{D^\bullet} \ ,
    \]
    and the desired equation follows from a simple computation. For general monomials $A,B$, the degree condition forces all $m$ homotopies to act on the pair of factors selected by the unique $d_{\mathrm{CE}}$. After cancellations, the only surviving terms are in the first summand of $\tb^{(m)}$, and where $d_\mathrm{CE}$ acts on a factor of $A$ and a factor of $B$, and where every internal homotopy acts on the same factor of the first input, transforming that form into $\alpha_{D^\bullet}$. The outgoing homotopy gives the same coefficient as on generators. Summing over these pairs with the appropriate sign gives precisely the biderivation extension $\{A,B\}_{(1-m)}$ of the Lie bracket.
\end{proof}

\noindent In dimension two this is the scalar of \cite[(4.44)]{BSV24}. We now show that in any dimension, the alternating sum of these scalars over a cross-polytope configuration has a geometric interpretation that fixes both its magnitude and its sign, making it independent of choices.
\begin{lemma}
    \label{lem:diff-geom}
    For a positively oriented cross-polytope configuration $\mathcal C$ as in Definition~\ref{def:cross-polytope-configuration}, with $T_\varepsilon=T(D^1_{\varepsilon_1},\ldots,D^m_{\varepsilon_m})$, one has
    \[
        \sum_{\varepsilon\in\{1,2\}^m}
          (-1)^{|\varepsilon|}\varphi_0(T_\varepsilon)=(-1)^m \ .
    \]
\end{lemma}
\begin{proof}
    For each tree $T_\varepsilon$, we denote $\alpha_\varepsilon \coloneqq k_{D^m_{\varepsilon_m}}\cdots k_{D^2_{\varepsilon_2}}(\omega_{D^1_{\varepsilon_1}}) \in \Omega^1_c(\R^m)$, and define recursively
    \[
        \eta_1\coloneqq\omega_{D^1_1}-\omega_{D^1_2}\in\Omega_c^m(\R^m)
    \]
    and, for $2\leqslant r\leqslant m$,
    \[
        \eta_r \coloneqq \left( k_{D^r_1}-k_{D^r_2}\right)\eta_{r-1} \ .
    \]
    We note that we have
    \[
        \eta_m =\sum_{\varepsilon\in\{1,2\}^m} (-1)^{|\varepsilon|}
        \alpha_\varepsilon \ ,
    \]
    and that $\eta_r$ is closed for $1\leqslant r\leqslant m$; this can be proven by induction from the fact that $\eta_1$ is closed and has integral zero. 
    Let us now denote by
    \[ 
        c = \sum_{\varepsilon\in\{1,2\}^m} 
        (-1)^{|\varepsilon|} \varphi^{(m)}_0(T_\varepsilon)
    \]
    the numerical value we are trying to compute. Using the definitions above we rewrite this as
    \[ 
        c = \frac12\int_{\R^m}
        \big(
        \phi
        (\omega_{\widetilde D}+\omega_{\R^m})
        +
        (-1)^m
        k_{\R^m}(\omega_{\widetilde D})
        \wedge d_{\mathrm{dR}}\phi
        \big) \ , 
    \]
    and using the homotopy $k_{\R^m}$ together with Stokes gives
    \[
        c = \int_{\R^m} \phi\, \omega_{\widetilde D} \ .
    \]
    We can further rewrite this integral by moving the form $\omega_{\widetilde D}$ away from the support of $\phi$. Write
    \[
        U_r\coloneqq D^r_1\cup D^r_2 \ ,\qquad 1\leqslant r\leqslant m \ ,
        \qquad U\coloneqq U_m \ ,
    \]
    so that $U_{r-1}\subset D^r_1\cap D^r_2$ and $\operatorname{supp}(\eta_r)\subset U_r$. Choose an open Euclidean ball $B$ containing $\overline D$. Since $d_{\mathrm{dR}}\phi=\eta_m$ has compact support contained in $U$, the function $\phi$ is constant on each component of $\R^m\setminus\operatorname{supp}(\eta_m)$, that is, it vanishes outside $B$ and on a collar of its boundary, and has a constant value on $\widetilde D$.     We now pick a form $\xi\in\Omega_c^{m-1}(\R^m)$ such that
    \[
        d_{\mathrm{dR}}\xi=\omega'-\omega_{\widetilde D} \ ,
    \]
    where $\omega'$ is supported outside $\overline B$. 
    Note that $\xi$ is closed on $X=D\setminus\overline{\widetilde D}$. We then get
    \[
        c=\int_{\R^m}\phi\,(\omega'-d_{\mathrm{dR}}\xi)
        =(-1)^{m-1} \int_{\R^m}\xi \wedge \eta_m \ ,
    \]
    by using Stokes and noting that $\operatorname{supp}(\omega')\cap\operatorname{supp}(\phi)=\varnothing$. We will compute the integral above by relating it to an integral over the boundary of the cross polytope that corresponds to the configuration $\mathcal{C}$.
    \medskip

    \noindent \textit{Geometric realization of the cross polytope:} Let us now realize the disk chains geometrically. Choose a point $p_V$ in each disk $V$ occurring in the configuration, and realize each simplex $\sigma=[V_0\subsetneq\cdots\subsetneq V_q]$ as a smooth singular chain $S(\sigma)$ with compact support in $V_q$, with vertices at the points $p_V$ and compatibly with the differential, that is:
    \[
        S[V_i]=p_{V_i} \ , \qquad \partial S(\sigma)=S(\partial\sigma) \ .
    \]
    These chains can be chosen by induction on dimension, using contractibility of each disk. We now use these chains to geometrically realize the cross polytope. We set $Z_1=p_{D^1_1}-p_{D^1_2}$ and
    \[
        B^r_j\coloneqq\sum_{\varepsilon\in\{1,2\}^{r-1}}(-1)^{|\varepsilon|}
        S\bigl[D^1_{\varepsilon_1}\subsetneq\cdots\subsetneq D^{r-1}_{\varepsilon_{r-1}}\subsetneq D^r_j\bigr] \ ,
        \qquad Z_r\coloneqq B^r_1-B^r_2 \ ,
    \]
    for all $2\leqslant r\leqslant m$ and $j \in \{1,2\}$. Each $B^r_j$ is an $(r-1)$-chain supported in $D^r_j$, and we calculate
    \[
        \partial B^r_j=(-1)^{r-1}Z_{r-1} \ .
    \]
    Geometrically, the first step gives two paths with the same endpoints, the next gives two surfaces filling the resulting cycle, and so on. In particular, the signed sum of top-dimensional chains is
    \[ 
        Z_m=S(z_{\mathcal C}) \ , 
    \]
    where $z_{\mathcal{C}} \in C_{m-1}(N_\bullet\mathscr{D}(D;\widetilde{D}))$ is the cycle of disk chains corresponding to the boundary of the cross polytope. Note that the positive orientation of $\mathcal C$ means that $Z_m$ represents the positive linking class around $\widetilde D$ in $X$.
    \medskip

    \noindent \textit{Iterated Stokes:} We claim that, for every closed $(r-1)$-form $\theta$ on the open $U_r=D^r_1\cup D^r_2$, we have
    \[
        \int_{U_r}\theta\wedge\eta_r=\int_{Z_r}\theta \ .
    \]
    We prove this by induction on $r$. For $r=1$, the function $\theta$ is constant on each of the disks $D^1_j$, so
    \[
        \int_{U_1}\theta\,(\omega_{D^1_1}-\omega_{D^1_2})
        =\theta(p_{D^1_1})-\theta(p_{D^1_2})=\int_{Z_1}\theta \ .
    \]
    For $r\geqslant2$, since each $D^r_j$ is contractible, we can choose $\lambda_j\in\Omega^{r-2}(D^r_j)$ such that $d_{\mathrm{dR}}\lambda_j=\theta|_{D^r_j}$, and using Stokes we get
    \[
        \begin{aligned}
            \int_{U_r}\theta\wedge\eta_r
            &=\int_{D^r_1}d_{\mathrm{dR}}\lambda_1\wedge k_{D^r_1}(\eta_{r-1})
            -\int_{D^r_2}d_{\mathrm{dR}}\lambda_2\wedge k_{D^r_2}(\eta_{r-1})\\
            &=(-1)^{r-1}\int_{U_{r-1}}(\lambda_1-\lambda_2)\wedge\eta_{r-1} \ .
        \end{aligned}
    \]
    from which we conclude by induction the desired identity. Applying this identity with $r=m$ and $\theta=\xi|_{U_m}$ to the expression for $c$ obtained by Stokes above gives
    \[
        c=(-1)^{m-1}\int_{Z_m}\xi = (-1)^{m-1} \int_E(\omega'-\omega_{\widetilde D})= (-1)^m \ ,
    \]
    where $E$ is a smooth disk satisfying $\widetilde D\Subset E\Subset D$. The equation above follows from the fact that $\xi$ is closed on $X$ and $Z_m$ represents the positive linking class and is homologous to the boundary of the smooth disk $E$.
\end{proof}

\begin{theorem}
    \label{thm:factorization-envelope-higher-bracket}
    Let $m\geqslant2$. The invariant associated with the factorization envelope $\mathfrak U\mathfrak g^{\R^m}$ is related to the canonical degree $1-m$ Poisson bracket on $\Sym(\mathfrak g[1-m])$ by
    \[
        L^{(m)}
        =
        (-1)^m\{-,-\}_{(1-m)} \ .
    \]
    Consequently, if $\mathfrak g$ is non-abelian, weight $m$ is the first non-zero higher weight of the transferred structure.
\end{theorem}

\begin{proof}
    Summing the formula for $\tb^{(m)}$ on linear trees over the cross-polytope cycle and applying \cref{lem:diff-geom} gives
    \[
        L^{(m)}
        =
        (-1)^m\{-,-\}_{(1-m)} \ .
    \]
    This is nonzero if $\mathfrak g$ is non-abelian, which implies $[\tb^{(m)}]\neq0$, by Proposition~\ref{prop:higher-poisson-invariant}. Together with Proposition~\ref{prop:factorization-envelope-weight-vanishing}, this proves the assertion about the first non-zero higher weight.
\end{proof}

\noindent For $m\geqslant3$, the first-order class of $\mathfrak U\mathfrak g^{\R^m}$ vanishes \cite[Section~4.3]{BSV24}, and non-trivial higher obstruction classes were expected \cite[Section~1]{BSV24}; \cref{thm:factorization-envelope-higher-bracket} confirms this at weight $m$ for non-abelian $\mathfrak g$. In conclusion, combining \cref{prop:abelian-formal} with \cref{thm:factorization-envelope-higher-bracket}, we conclude:
\begin{corollary}
    \label{cor:factorization-envelope-equivalence}
Let $m\geqslant2$. For the factorization envelope $\mathfrak U\mathfrak g^{\R^m}$, the following are equivalent:
    \begin{enumerate}
        \item $\mathfrak{g}$ is abelian,
        \item $L^{(m)}=0$,
        \item $\mathfrak U\mathfrak g^{\R^m}$ is formal.
    \end{enumerate}
\end{corollary}

\subsection{Factorization envelopes of $L_\infty$-algebras}
\label{subsec:dg-Linfty-envelopes}
We now generalize the discussion of the last subsection to the setting where $\mathfrak g$ denotes an arbitrary (uncurved) $L_\infty$-algebra with operations $\ell_k$ of degree $2-k$, and differential $\ell_1$. This subsumes the case examined above as well as any dg Lie algebra, seen as a $L_\infty$-algebra with trivial higher operations $\ell_k=0$ for $k\geqslant3$. Similarly to the classical case, the envelope is
\begin{equation}
    \label{eq:general-factorization-envelope}
    \mathfrak U^{(m)}_{\mathfrak g}(D)
    =\left(\Sym\bigl((\mathfrak g\otimes\Omega_c^\bullet(D))[1]\bigr),
    d_{\mathrm{lin}}+\delta\right),
    \qquad \delta=\sum_{k\geqslant2}d_k \ .
\end{equation}
Here $d_{\mathrm{lin}}$ is induced by the shifted tensor differential, and $d_k$ is the Chevalley--Eilenberg term replacing $k$ generators by their $\ell_k$-operation. Another difference is that one needs to keep track of internal degrees for sign conventions: for homogeneous $a,b\in\mathfrak g$, the binary differential is
\[
    d_2\bigl((a\otimes\omega)(b\otimes\eta)\bigr)
    =(-1)^{|a|+|\omega|_{\mathrm{dR}}(|b|+1)}
    \ell_2(a,b)\otimes\omega\wedge\eta \ ,
\]
combining the shifted Chevalley--Eilenberg sign with the Koszul sign for $b$ and $\omega$. Choose a strong deformation retract $(i_{\mathfrak g},p_{\mathfrak g},h_{\mathfrak g})$ of ${\mathfrak g}$ to its cohomology, using the sign convention
\[
\ell_1 h_{\mathfrak{g}} + h_{\mathfrak{g}} \ell_1
= \id - i_{\mathfrak{g}} p_{\mathfrak{g}} \ ,
\]
and extend the linear retract
\[
\begin{aligned}
i_D(a)&=i_{\mathfrak g}(a)\otimes\omega_D \ ,\\
p_D(a\otimes\omega)&=p_{\mathfrak g}(a)\int_D\omega \ ,\\
h_D(a\otimes\omega)&=(-1)^{|a|}a\otimes k_D(\omega)
+h_{\mathfrak g}(a)\otimes\omega_D\int_D\omega \ .
\end{aligned}
\]
to the unperturbed complex $\bigl(\Sym((\mathfrak g\otimes\Omega_c^\bullet(D))[1]),d_{\mathrm{lin}}\bigr)$ by the same symmetrized tensor trick. As before, the integrals vanish outside form degree $m$. Although $h_{\mathfrak g}$ follows the opposite sign convention, the shift ensures that the resulting homotopy satisfies
\[
d_{\mathrm{lin}}h_D+h_Dd_{\mathrm{lin}}=i_Dp_D-\id \ ,
\]
as required in Section~\ref{sec:obstruction-sequences}, together with the side conditions.
\medskip

\noindent Since $d_k$ lowers symmetric degree by $k-1$, the homological perturbation lemma applies to $\delta$, and gives the retract by the same formula with $d_{\mathrm{CE}}$ replaced by $\delta$. As in the previous subsection, $\delta i_D=0$, since every term of $\delta$ wedges at least two top forms on $D$. Thus the inclusion is unchanged, the transferred differential is zero, and on every disk $D$,
\[
H\mathfrak U^{(m)}_{\mathfrak g}(D)
\cong\Sym\bigl(H^\bullet(\mathfrak g)[1-m]\bigr) \ .
\]
From now on, let $m\geqslant2$. The argument of \cite[Section~4.2]{BSV24} then shows that the underlying strict prefactorization algebra structure is the graded commutative product on $\Sym\bigl(H^\bullet(\mathfrak g)[1-m]\bigr)$. The following theorem generalizes the results of the previous subsection. Its proof is the same, except that elements now have many degrees and the retract of $\mathfrak g$ contributes additional homotopies.
\begin{theorem}
\label{thm:coefficient-envelope-first-invariant}
For $m\geqslant2$, the original and improved transfers of \eqref{eq:general-factorization-envelope} satisfy
\[
\widehat\beta^{(2)}=\cdots=\widehat\beta^{(m-1)}=0 \ ,
\qquad
\widetilde\beta^{(2)}=\cdots=\widetilde\beta^{(m-1)}=0 \ ,
\]
and the associated higher Poisson invariant is
\[
L^{(m)}=(-1)^m\{-,-\}_{(1-m)} \ ,
\]
where the degree $1-m$ bracket on $\Sym(H^\bullet(\mathfrak g)[1-m])$ is extended from
\[
\{a,b\}_{(1-m)}
=(-1)^{(m-1)|a|}\ell_2^H(a,b) \ ,\qquad a,b\in H^\bullet(\mathfrak g) \ ,
\]
as a biderivation. Here $\ell_2^H$ is the induced Lie bracket on $H^\bullet(\mathfrak g)$, given by
\[
\ell_2^H(a,b)=p_{\mathfrak g}\ell_2(i_{\mathfrak g}a,i_{\mathfrak g}b) \ .
\]
Therefore, if $\ell_2^H\neq0$, weight $m$ is the first non-zero higher weight, and $[\widetilde\beta^{(m)}]\neq0$.
\end{theorem}

\begin{proof}
    The first part of this proof uses an alternative degree-counting argument from the proof of \cref{prop:factorization-envelope-weight-vanishing}, to account for all the terms in the Chevalley--Eilenberg differential as well as the contracting homotopy of $\mathfrak{g}$.
    \medskip
    
    \noindent Consider a summand in the perturbative expansion for $\widehat{\beta}^{(r)}$. Let $q_k$ count occurrences of $d_k$, and let $u$ and $v$ count $h_{\mathfrak{g}}$ and de Rham homotopies (that is, $k_E$ for some disk $E$), respectively. These homotopies appear once for each internal edge and once each time one applies a $d_k$, so
        \[
        u+v=r-1+\sum_{k\geqslant2}q_k \ .
        \]
    Each $d_k$ preserves total base form degree and removes $k-1$ factors, whereas only the de Rham homotopies lower that degree. Since the input is a product of $m$-forms, and the final projection only survives if all remaining factors have base form degree $m$, the number of de Rham homotopies must satisfy
        \[
        v=m\sum_{k\geqslant2}(k-1)q_k \ ,
        \]
    giving
        \[
        r-1=\sum_{k\geqslant2}\bigl((k-1)m-1\bigr)q_k+u \ .
        \]
    If $r<m$, every $q_k$ vanishes, so all homotopies are coefficient homotopies. For $r\geqslant2$, the one on a lowest internal edge is applied to the image of $i_{\mathfrak g}$, and the term vanishes since $h_{\mathfrak g}i_{\mathfrak g}=0$. Proposition~\ref{prop:higher-weight-comparison} gives the same vanishing for $\tb$.
    \medskip
    
    \noindent The rest of the proof follows in the same fashion as in the previous subsection, with the only difference being the signs one gets when applying a homotopy. For homogeneous $a\in H^\bullet(\mathfrak g)$, put $A=i_{\mathfrak g}(a)$. Then
        \[
        h_{D^m}\cdots h_{D^2}(A\otimes\omega_{D^1})
        =(-1)^{(m-1)|a|}A\otimes\alpha_{D^\bullet} \ .
        \]
    This accounts for the overall sign in the shifted coefficient bracket. The only differential-geometric calculation one needs to do is the same as before, that is, \cref{lem:diff-geom}.
\end{proof}

\begin{corollary}
    \label{cor:coefficient-envelope-weight-m-vanishing}
    Let $m\geqslant2$ and consider the factorization envelope \eqref{eq:general-factorization-envelope} of an uncurved $L_\infty$-algebra $\mathfrak g$, with the improved transfer constructed above. If $L^{(m)}=0$, equivalently $\ell_2^H=0$, then $\tb^{(m)}=0$.
\end{corollary}
\begin{proof}
    By the weight-counting identity in the proof of the preceding theorem, every surviving term of $\widehat\beta^{(m)}$ contains exactly one $d_2$ and no coefficient homotopy. For arbitrary monomial inputs, its coefficient part therefore contains
    \[
        p_{\mathfrak g}\ell_2(i_{\mathfrak g}a,i_{\mathfrak g}b)
        =\ell_2^H(a,b)
    \]
    for the two selected generator factors, while every other coefficient passes unchanged to cohomology. Thus $\ell_2^H=0$ implies $\widehat\beta^{(m)}=0$, and Proposition~\ref{prop:higher-weight-comparison} gives $\tb^{(m)}=0$. The equivalence with $L^{(m)}=0$ follows from Theorem~\ref{thm:coefficient-envelope-first-invariant}.
\end{proof}

\begin{remark}
    When $\ell_2^H=0$, Corollary~\ref{cor:coefficient-envelope-weight-m-vanishing} gives $\tb^{(m)}=0$, but this does not imply triviality of the components $\tb^{(r)}$ for $r>m$, which may depend on higher $L_\infty$ operations. This is a fundamental difference between this case and the case where $\mathfrak{g}$ is an ordinary Lie algebra. When there are higher structures, we do not get a direct generalization of the `if and only if' statements of \cref{cor:factorization-envelope-equivalence}.
\end{remark}

\subsection{Compactified factorization envelopes}
\label{subsec:compactified-factorization-envelopes}
We would like to specialize the previous subsection to a case of particular interest in physics, where the dg Lie algebra is given by
\[ 
    \mathfrak g\otimes\Omega^\bullet(N), \quad 
    [X\otimes\alpha,Y\otimes\beta] =[X,Y]\otimes\alpha\wedge\beta \ , 
\]
where $N$ is a closed manifold and $\mathfrak g$ is a Lie algebra of finite dimension concentrated in degree zero. The factorization envelope \eqref{eq:general-factorization-envelope} of this dgla can be seen as a polynomial model for the pushdown of the factorization envelope of $\mathfrak{g}$ on the total space $\R^m\times N$ of a trivial fiber bundle, in other words, the theory compactified over $N$. Explicitly, it is given by:
\[
    \mathfrak U^{(m)}_{\mathfrak g,N}(D)
    =\left(
    \Sym\bigl((\mathfrak g\otimes\Omega^\bullet(N)\otimes\Omega_c^\bullet(D))[1]\bigr),
    d_{\mathrm{lin}}+d_{\mathrm{CE}}
    \right).
\]
Choose a strong deformation retract $(i_N,p_N,h_N)$ of $\Omega^\bullet(N)$ onto its cohomology, with $d_Nh_N+h_Nd_N=\id-i_Np_N$, and put $V_N=\mathfrak g\otimes H^\bullet_{\mathrm{dR}}(N)$, endowed with the Lie bracket
\[
    [X\otimes[\alpha],Y\otimes[\beta]]_{\cup}
    =[X,Y]\otimes[\alpha\wedge\beta] \ .
\]
Then
\[
    H\mathfrak U^{(m)}_{\mathfrak g,N}
    \cong\Sym\left(V_N[1-m]\right) \ .
\]
Spelling out Theorem~\ref{thm:coefficient-envelope-first-invariant}, we have vanishing of $\tb^{(r)}$ for $2\leqslant r<m$ and
\[
    L^{(m)}=(-1)^m\{-,-\}_{(1-m)} \ ,
\]
where the biderivation on the right is determined by the cup product of cohomology classes of $N$ by
\[
    \{X\otimes[\alpha],Y\otimes[\beta]\}_{(1-m)}
    =(-1)^{(m-1)|\alpha|}[X,Y]\otimes[\alpha\wedge\beta] \ .
\]
We note that, upon fixing $\mathfrak g$ and the `de Rham homotopies' $k_D$, along the base $\R^m$, $\tb^{(m)}$ only depends on the fiber $N$ through the bracket above. This follows from the argument in the proof of \cref{cor:coefficient-envelope-weight-m-vanishing}. Note that, for this argument characterizing the dependence of $L^{(m)}$ on $N$, no formality assumption on this manifold is ever needed, and its Massey products never appear. Moreover, the equivalence of \cref{cor:factorization-envelope-equivalence} extends to compactifications:
\begin{proposition}
    \label{prop:compactified-abelian-formal}
    Let $m\geqslant2$ and let $N$ be nonempty. Then the following are equivalent:
    \begin{enumerate}
        \item $\mathfrak{g}$ is abelian,
        \item $L^{(m)}=0$,
        \item the compactified factorization envelope $\mathfrak U^{(m)}_{\mathfrak g,N}$ is formal.
    \end{enumerate}
    Moreover, if $\mathfrak{g}$ is abelian, then $\widehat\beta^{(r)} = \tb^{(r)} = 0$ for every $r\geqslant2$.
\end{proposition}
\begin{proof}
    If $\mathfrak g$ is abelian, the bracket of $\mathfrak g\otimes\Omega^\bullet(N)$ vanishes, so every $q_k$ is zero in the weight count of the proof of \cref{thm:coefficient-envelope-first-invariant}, and the argument given there for $r<m$ applies in every weight $r\geqslant2$. This proves the last assertion, and hence formality. Formality implies $L^{(m)}=0$ by \cref{prop:higher-poisson-invariant}. Finally, if $[X,Y]\neq0$ for some $X,Y\in\mathfrak g$, then $\{X\otimes[1],Y\otimes[1]\}_{(1-m)}=[X,Y]\otimes[1]$ is nonzero, since $[1]\neq0$ in $H^0_{\mathrm{dR}}(N)$, so $L^{(m)}\neq0$ by \cref{thm:coefficient-envelope-first-invariant}.
\end{proof}
\begin{remark}
    For non-abelian $\mathfrak g$, the envelope $\mathfrak U^{(m)}_{\mathfrak g,N}$ is never formal, and this is already detected by $L^{(m)}$, independently of $N$. Nevertheless, higher-weight components may depend on the higher structure of $N$. In the weight count of the proof of \cref{thm:coefficient-envelope-first-invariant}, a coefficient homotopy $h_N$ can only survive in terms with $q_2\geqslant2$. 
    The first such terms have weight $2m$, with $q_2=2$ and $u=1$, and involve the expressions $p_N\bigl(h_N(\alpha\wedge\beta)\wedge\gamma\bigr)$ for closed forms $\alpha,\beta,\gamma$, from which the triple Massey products of $N$ are built. For instance, the Heisenberg nilmanifold has a nontrivial triple Massey product \cite[Section~2]{FM03}; its minimal model is the algebra $A$ of \cref{prop:massey-counterexample}; we conjecture that this manifold, together with a non-abelian Lie algebra, such as $\mathfrak{sl}_2$, gives an example in which this Massey product is detected at weight $2m$.
\end{remark}

\section{\textcolor{bordeau}{Symmetric pairings}}
\label{sec:complexes-with-pairing}
\noindent We now consider a second type of examples, very similar to the factorization envelope construction that we studied in the previous section, but where the perturbation to the differential comes from a \emph{pairing}, that is, an operation that contracts two generators to a scalar (instead of a linear term, as happened in that previous type of example). Special cases of this construction give prefactorization algebras associated to free BF-type theories, abelian Chern--Simons theories and their compactifications. We adjoin a formal parameter $\hbar$ of cohomological degree zero for a closer connection to the physics literature. In that way, observables and all transferred operations are defined over $\R[\hbar]$.

\subsection{Prefactorization algebras from symmetric pairings}
Fix $m\geqslant1$, and let $V$ be a cochain complex equipped with a closed graded symmetric pairing
\[
    B_V\colon V\otimes V\longrightarrow\R[m+1] \ ,
\]
which is allowed to be degenerate. That is, we have a scalar pairing satisfying
\[
    B_V(v,w)=(-1)^{|v||w|}B_V(w,v) \ ,
    \qquad
    B_V(d_Vv,w)+(-1)^{|v|}B_V(v,d_Vw)=0 \ ,
\]
and we note that a nonzero pairing requires $|v|+|w|=-m-1$.

For each disk $D$ we look at the space of $V$-valued forms $\mathcal E_V(D)=V\otimes\Omega_c^\bullet(D)$ endowed with the `linear' differential
\[
    d_{\mathrm{lin}}(v\otimes\alpha)
    =d_Vv\otimes\alpha+(-1)^{|v|}v\otimes d_{\mathrm{dR}}\alpha
\]
and a pairing induced from $B_V$
\[
    B_D(v\otimes\alpha,w\otimes\beta)
    =(-1)^{|\alpha||w|}B_V(v,w)\int_D\alpha\wedge\beta.
\]
This is a closed graded symmetric pairing of degree one, which we extend to
\[
    \Delta_V(\xi_1\cdots\xi_a)
    =\sum_{i<j}(-1)^{\epsilon_{ij}}B_D(\xi_i,\xi_j)
    \prod_{k\neq i,j}\xi_k \ , 
\]
where $(-1)^{\epsilon_{ij}}$ is the Koszul sign for bringing $\xi_i,\xi_j$, in that order, to the front. This second-order operator vanishes on constants and linear generators, and satisfies
\[
    \Delta_V^2=0,\qquad
    d_{\mathrm{lin}}\Delta_V+\Delta_Vd_{\mathrm{lin}}=0 \ .
\]
Consequently, for each disk $D$,
\[
    \mathfrak F_V(D)
    \coloneqq\left(\Sym(\mathcal E_V(D))[\hbar],d_{\mathrm{lin}}+\hbar\Delta_V\right)
\]
is a cochain complex. We interpret its specialization at $\hbar=0$ as the classical complex of observables on $D$, and $\hbar\Delta_V$ as the BV quantum correction to the differential. 

\noindent As in the factorization envelope examples, we can make the complex above into a prefactorization algebra by using extension by zero and symmetric multiplication, over $\R[\hbar]$. Note that these maps remain cochain maps for the perturbed differential: the second-order operator $\Delta_V$ would only ever fail to be a derivation by terms pairing two forms, but since this is always applied to forms with disjoint supports, that failure vanishes.
\medskip

\noindent Let $W=H^\bullet(V)$
and choose a retract $(i_V,p_V,h_V)$ from $V$ to $W$, with
\[
    d_Vh_V+h_Vd_V=i_Vp_V-\id \ ,
\]
as in Section~\ref{sec:obstruction-sequences}, and de Rham homotopies $k_D$ as in \cref{subsec:factorization-envelopes}. On the linear complexes take
\begin{align*}
    i_D(v)&=i_V(v)\otimes\omega_D \ ,\\
    p_D(v\otimes\omega)&= p_V(v) \int_D\omega \ ,\\
    h_D(v\otimes\omega)&=(-1)^{|v|+1}v\otimes k_D(\omega)
    +h_V(v)\otimes\omega_D\int_D\omega \ ,
\end{align*}
so that we have, as desired,
\[
    d_{\mathrm{lin}}h_D+h_Dd_{\mathrm{lin}}=i_Dp_D-\id \ .
\]
We extend these maps to symmetric powers using the same tensor trick and perturb by the quantum correction $\hbar\Delta_V$. The resulting maps $\widehat p_D,\widehat h_D$ are given by the same expansion as before, with $d_{\mathrm{CE}}$ replaced by $\hbar\Delta_V$. The inclusion is unchanged and the transferred differential is zero: $\Delta_Vi_D=0$, since two top forms on $D$ wedge to zero. Thus,
\[
    H\mathfrak F_V(D)\cong\Sym(W[-m])[\hbar]
\]
and the prefactorization algebra $\mathfrak F_V$ is locally constant.

\medskip

\noindent In what concerns its transfers $\widehat{\beta}$ and $\tb$, as in the factorization envelope examples, dimension $m=1$ differs from the $m\geqslant2$ case. We first focus on the latter and return to dimension one below.
\begin{theorem}
    \label{thm:paired-complex-first-invariant}
    For $m\geqslant2$, the original and improved transfer components $\widehat\beta^{(r)}$ and $\widetilde\beta^{(r)}$ can be nonzero only if
    \[
        r=1+q(m-1) \ ,\qquad q\in\mathbb Z_{\geqslant0} \ .
    \]
    The strict product $\widetilde\beta^{(1)}$ is the usual product on $\Sym(W[-m])[\hbar]$, and the associated higher Poisson invariant is the $\R[\hbar]$-bilinear degree $1-m$ biderivation on $\Sym(W[-m])[\hbar]$ determined by its value on linear monomials
    \[
        L^{(m)}(v, w)=\hbar(-1)^{m|v|}B_W(v,w) \ ,
    \]
    where $B_W(v,w)\coloneqq B_V\big(i_V(v),i_V(w)\big)$ is the induced pairing on $W$, which is independent of the choice of $i_V$ since $B_V$ is closed. If $B_W\neq0$, the weight-$m$ obstruction class is nonzero, and this is the first nonzero higher weight.
\end{theorem}
\begin{proof}
    The entire proof follows in the same fashion as in the preceding examples, so here we only note the differences here. Firstly, coefficient labels begin in $\operatorname{im}(i_V)$. Products retain them as separate factors and contractions remove pairs, so no operation creates a new coefficient label. Thus every term containing an occurrence of $h_V$ vanishes.
    \medskip
    
    \noindent The weight restriction follows from a similar degree count as in \cref{prop:factorization-envelope-weight-vanishing}, and
    for $r=1$ there are no contractions, giving the usual product. As for the sign of the invariant, the linear-tree calculation gives
    \[
        \tb^{(m)}(T)(v,w)
        =\hbar(-1)^{m+m|v|}\varphi_0^{(m)}(T)B_W(v,w) \ ,
    \]
    so the sum over a positively oriented cross-polytope configuration reduces again to the calculation of \cref{lem:diff-geom}, giving the stated formula. We extend this to general monomial inputs by the biderivation property of $L^{(m)}$, from Theorem~\ref{thm:higher-poisson-bracket}.
\end{proof}

\noindent We note a fundamental difference between this case and the examples in the previous section: as argued in the proof above, any component of $\widetilde\beta$ is independent of the choice of internal homotopy $h_V$. As a result, none of them can depend on the higher structure on $W$. Another way of seeing this is by using the following result.
\begin{proposition}
    \label{prop:paired-coefficient-reduction}
    Equip $W$ with its induced pairing $B_W$ and zero differential. Any choice of closed representatives $i_V\colon W \to V$ induces an objectwise quasi-isomorphism of prefactorization algebras $\mathfrak F_W\xlongrightarrow{\simeq}\mathfrak F_V$.
\end{proposition}
\begin{proof}
    We have $B_V(i_Vv,i_Vw)=B_W(v,w)$, so the $\R[\hbar]$-linear extension of $\Sym(i_V\otimes\id)$ commutes with the linear differentials, contractions, and prefactorization maps, giving a morphism of prefactorization algebras
    \[ 
        \mathfrak F_{W}\to\mathfrak F_V \ ,
    \]
    and all that remains is to check the quasi-isomorphism. This can be done by filtering by symmetric degree and checking on each graded piece, since the contraction lowers this degree by two and the map on each associated graded piece is a finite symmetric power of quasi-isomorphisms.
\end{proof}

\noindent In other words, the weak equivalence type of our construction depends only on $H^\bullet(V)$ and its pairing, even beyond this first invariant $L^{(m)}$. As a result, if the induced pairing on $W$ is zero, then $\mathfrak F_V$ is a formal prefactorization algebra.
\vspace{0.5cm}

\noindent In the remaining case $m=1$, the pairing enters the product itself; in analogy with the factorization envelope case, where we obtain the Gutt star-product deforming the symmetric algebra to $\mathfrak{Ug}$, here we get a Moyal star product.

\begin{proposition}
    \label{prop:paired-complex-moyal}
    In dimension $m=1$, the higher components $\tb^{(r)}$ vanish for all $r\geqslant2$, and $\mathfrak F_V$ is formal. For $D_1$ to the left of $D_2$, the strict product on $H$ is the graded Moyal product
    \[
        f\star_{\mathrm M}g
         =\mathrm{mult}\circ\exp\!\left(\frac{\hbar}{2}\partial_\Pi\right)(f\otimes g) \ ,
        \qquad
        \Pi(\bar v,\bar w)=(-1)^{|v|+1}B_W(v,w) \ .
    \]
    Here $\mathrm{mult}$ denotes symmetric multiplication, and $\partial_\Pi$ sums over contractions of one generator of $f$ with one generator of $g$ using $\Pi$ and   Koszul signs. The commutator of two generators is
    \[
       \bar v\star_{\mathrm M}\bar w
       -(-1)^{(|v|+1)(|w|+1)}\bar w\star_{\mathrm M}\bar v
       =\hbar\,\Pi(\bar v,\bar w) \ .
    \]
    Thus this associative product quantizes the constant Poisson bracket determined by $\Pi$.
\end{proposition}

\begin{proof}
    The product map $\mu_{(D_1,D_2)}$ is given by
    \[ \mu_{(D_1,D_2)}(f,g) = \widehat{\beta}\bigl(\op{\R}{D_1,D_2}\bigr)(f,g) = \widehat p_{\R}((f\otimes \omega_1)(f\otimes \omega_2)) \ , \]
    where $\omega_j=\omega_{D_j}$. Let $C$ be the second-order operator on $\Sym(\mathcal E_V(\R))[\hbar]$ given on quadratic monomials by $C(\xi\eta)=\Delta_Vh_{\R}(\xi\eta)$, and extended to higher monomials by summing over pairs, as in the definition of $\Delta_V$. Thus we have an expression
    \[ \widehat p_{\R} = \sum_{j\geqslant 0}\hbar^j p_\R C^j \]
    for the perturbed projection. Writing $F_j(t)=k_\R(\omega_j) = \int_{-\infty}^t\omega_j$, we then calculate, using the symmetrized homotopy and an integration by parts, that
    \[
        C\bigl((v\otimes\omega_1)(w\otimes\omega_2)\bigr)
        =\frac12(-1)^{|v|+1}B_W(v,w)\int_\R\bigl(F_1\omega_2-F_2\omega_1\bigr)
        =\frac12\,\Pi(\bar v,\bar w) \ ,
    \]
    since $F_1=1$ on $D_2$ and $F_2=0$ on $D_1$. The same formula shows that $C$ vanishes on two generators carrying the same form. Note that the last map $p_\R$ integrates the remaining top forms to one; the desired expression for the product then follows from matching the symmetry coefficients for the sum $\sum_{j\geqslant 0} C^j$ acting on symmetric powers, and for $\exp(\hbar\,\partial_\Pi/2)$ acting on the tensor.
\end{proof}

\noindent Let us present some concrete examples. Taking a symplectic vector space $(\mathfrak a,\kappa)$ gives a complex $V=\mathfrak a[1]$, with a symmetric pairing $B_V=\kappa$ of degree $m+1=2$. In this case, 
\[ H=\Sym(\mathfrak a)[\hbar]\quad \text{and} \quad \Pi=\kappa \ , \]
so $(H,\star)$ is the Weyl algebra of $(\mathfrak a,\kappa)$, recovering the complex of observables of topological quantum mechanics \cite[Section~5.1]{GRW21}.
\medskip

\noindent Similarly, for a closed oriented surface $\Sigma$, the complex $V=\Omega^\bullet(\Sigma)[2]$ endowed with the pairing
\[ B_V(\alpha_2,\beta_2)=\int_\Sigma\alpha\wedge\beta \]
models abelian Chern--Simons theory compactified to a line, and the proposition gives the Weyl algebra of $H^\bullet(\Sigma)[1]$ with the pairing induced by Poincar\'e duality, as in \cite[Section~5.3.3]{GRW21}.

\subsection{Free BF theory}
\label{subsec:free-BF}
One class of examples of the construction for complexes with pairing we study in this section comes from BF-type theories \cite{Hor89,CR01}. More precisely, we consider generalizations of the abelian BF theory with shift $k$ of \cite[Section~3]{CMR18}, which has fields in
\[ 
    \Omega^\bullet(M)[n-k-1] \oplus \Omega^\bullet(M)[k] \ , 
\]
where $M$ is some manifold of dimension $n$. We will consider a polynomial model of this theory on $M = \R^m \times N$, compactified on a closed smooth manifold $N$ of dimension $d=n-m$. Let us choose integers $s,t$ with $s+t=n+1$, and pick a finite-dimensional vector space $\mathfrak a$. We consider the coefficient complex
\[
    V=(\mathfrak a\otimes\Omega^\bullet(N))[s]
    \oplus(\mathfrak a^\vee\otimes\Omega^\bullet(N))[t] \ ,
\]
with pairing equal to zero on each individual summand, and otherwise determined by
\[
    B_V(x\otimes\alpha_s,\lambda\otimes\beta_t)
    =(-1)^{t|\alpha|}\lambda(x)\int_N\alpha\wedge\beta \ .
\]
Note that this pairing is of degree $s+t-d=m+1$, so we can consider the prefactorization algebra $\mathfrak F_V$. The cohomology of this prefactorization algebra is then given by
\[
    H\mathfrak F_{\mathrm{BF}}^{N;s,t}
    \cong\Sym\bigl(\mathfrak a\otimes H^\bullet_{\mathrm{dR}}(N)[s-m]
    \oplus\mathfrak a^\vee\otimes H^\bullet_{\mathrm{dR}}(N)[t-m]\bigr)[\hbar] \ .
\]
For $m\geqslant2$, Theorem~\ref{thm:paired-complex-first-invariant} gives its higher Poisson invariant, whose value on linear monomials is
\[
L^{(m)}(x\otimes[\alpha],\lambda\otimes[\beta])
=\hbar(-1)^{(t-m)|\alpha|+ m s}\lambda(x)\int_N\alpha\wedge\beta \ ,
\]
and zero on pairs from the same summand. Therefore, for $\mathfrak a\neq0$ we get a nonzero first invariant in every base dimension $m\geqslant2$.
\medskip

\noindent If instead $m=1$, so that $\dim N=n-1$, \cref{prop:paired-complex-moyal} applies, and we get a formal prefactorization algebra with the Moyal-deformed strict product. Since $(s-1)+(t-1)=\dim N$, Poincar\'e duality identifies $\mathfrak a^\vee\otimes H^\bullet_{\mathrm{dR}}(N)[t-1]$ with the dual of
\[ U=\mathfrak a\otimes H^\bullet_{\mathrm{dR}}(N)[s-1] \ , \]
and $\Pi$ with the canonical pairing on this space, up to sign. We relate this to the $k$-shifted abelian BF theory of \cite{CMR18} on $M = N \times \R$ by taking $s=n-k$ and $t=k+1$; then, $H\mathfrak F_V$ is the Weyl algebra of $U\oplus U^\vee$, the quantization of the cotangent space $T^*U$ associated to a `spatial slice' $N$.
\medskip

\noindent For the usual BF theory, we take $k=1$, so the degree-zero fields on $M = N\times \R$ become a connection one-form and an $(n-2)$-form. The resulting algebra is precisely the Weyl algebra of $T^*H^1(N;\mathfrak a)$, the reduced phase space of abelian BF theory on $M$ \cite{Hor89}.

\subsection{Higher abelian Chern--Simons theory and its compactifications}
\label{subsec:abelian-CS-examples}
Lastly, we look at special cases of the construction for complexes with symmetric pairing that only make sense for certain dimensions. We choose
\[ 
    n=4\ell+3 \ , \quad s=(n+1)/2=2\ell+2 \ ,
\]
and let $(\mathfrak a,\kappa)$ be a finite-dimensional vector space with a nondegenerate symmetric pairing. 
\medskip

\noindent Higher abelian Chern--Simons theory \cite{GRW21} has fields given by $\mathfrak a\otimes\Omega^\bullet(-)[s-1]$, which, compactified on a nonempty closed oriented manifold $N$ of dimension $d=n-m$, gives the coefficient complex with pairing
\[
    V=(\mathfrak a\otimes\Omega^\bullet(N))[s] \ ,
    \qquad
    B_V(x\otimes\alpha,y\otimes\beta)
    =\kappa(x,y)\int_N\alpha\wedge\beta \ .
\]
Since $s$ is even and $2s-d=m+1$, this is a pairing of the required degree and symmetry. The resulting prefactorization algebra $\mathfrak F_{\mathrm{CS}}^{N,\mathfrak a}$ is a polynomial model for the pushforward of higher abelian Chern--Simons theory.
\medskip

\noindent Theorem~\ref{thm:paired-complex-first-invariant} then gives
\[
    H\mathfrak F_{\mathrm{CS}}^{N,\mathfrak a}
    \cong\Sym\bigl(\mathfrak a\otimes H^\bullet_{\mathrm{dR}}(N)[s-m]\bigr)[\hbar] \ ,
\]
with the shifted Poincar\'e pairing
\[
    L^{(m)}(x\otimes[\alpha],y\otimes[\beta])
    =\hbar(-1)^{(s-m)|\alpha|}\kappa(x,y)\int_N\alpha\wedge\beta \ .
\]
\medskip

\noindent We now compare this situation with that analyzed in \cite[Section~5]{BSV24}. For $n=3$, $m=2$, $N=S^1$, and $\mathfrak a=\R$ with its standard pairing, this gives $L^{(2)}([1],[dt])=\hbar$, with $\int_{S^1}dt=1$, recovering their Theorem~5.1 at $\hbar=1$. When we take the trivial compactification over $N=\mathrm{pt}$, one gets at cohomology level the constant family
\[ 
    \Sym(\mathfrak a[-(2\ell+1)])[\hbar] \ ,
\]
with first invariant $L^{(n)}(x,y)=\hbar\kappa(x,y)$. For $n=3$, this is abelian Chern--Simons theory on $\R^3$, whose first-order class vanishes. As expected in that reference, there are indeed non-trivial higher obstruction classes, namely the weight-three class $L^{(3)}=\hbar\kappa \neq 0$.

\bibliographystyle{alpha} \bibliography{bib}

\end{document}